\documentclass[12pt,leqno,ngerman,english]{article}
\usepackage[T1]{fontenc}
\usepackage[latin9]{inputenc}
\usepackage[letterpaper]{geometry}
\usepackage{babel}
\usepackage{float}
\usepackage{textcomp}
\usepackage{mathtools}
\usepackage{amsmath}
\usepackage{amsthm}
\usepackage{amssymb}
\usepackage{cancel}
\usepackage{esint}
\usepackage[unicode=true,pdfusetitle,
 bookmarks=true,bookmarksnumbered=false,bookmarksopen=false,
 breaklinks=false,pdfborder={0 0 0},pdfborderstyle={},backref=false,colorlinks=false]
 {hyperref}
\usepackage{breakurl}

\makeatletter

\newcommand{\lyxmathsym}[1]{\ifmmode\begingroup\def\b@ld{bold}
  \text{\ifx\math@version\b@ld\bfseries\fi#1}\endgroup\else#1\fi}

\usepackage{titling}
\posttitle{\par\end{center}\vspace{0em}}

\renewcommand\paragraph{\@startsection{paragraph}{4}{\z@}%
   {-3.25ex\@plus -1ex \@minus -.2ex}%
   {1.5ex \@plus .2ex}%
   {\normalfont\normalsize\bfseries}}

\usepackage{array}
\usepackage[indentafter]{titlesec}
\titleformat{\section}{\normalfont\Large\bfseries}{\S\  \thesection.}{1em}{}

\date{}

\makeatother

\begin{document}

\title{On Riemann's Nachlass for Analytic Number Theory}

\author{Carl Ludwig Siegel (1932)\\
\\
Translated by Eric Barkan and David Sklar\thanks{Eric Barkan (ebarkan@jps.net), 98 Gazania Court, Novato, CA 94945,
and David Sklar (dsklar@sfsu.edu), 326 6th Avenue, San Francisco CA
94118}\vspace{-20mm}
}
\maketitle
\begin{abstract}
In 1859 Riemann (1826$\:$--$\,$1866) published his only paper on
number theory. In this eight-page paper he obtained a formula for
the number of primes less than or equal to a real number $x$, and
revealed the deep connection between the distribution of primes and
the zeros of an analytic function now called the Riemann Zeta Function.
In the early 1930\textquoteright s two related unpublished results
from 1859 were found in Riemann\textquoteright s very rough notes
by the great twentieth century mathematician and scholar, Carl Ludwig
Siegel (1896 \nolinebreak -- \nolinebreak 1981).

In his 1932 paper \foreignlanguage{ngerman}{\emph{Über Riemanns Nachlaß
zur analytischen Zahlentheorie}} (\emph{On Riemann\textquoteright s
Nachlass for Analytic Number Theory}) Siegel presents these unpublished
results and gives derivations he found in Riemann\textquoteright s
notes. The first is an asymptotic development, now called the Riemann-Siegel
formula, for efficiently computing values of the Riemann zeta function.
The second is a new integral representation of the zeta function.
These had not been rediscovered seventy years after Riemann. Thus,
in 1932 the importance of Siegel\textquoteright s paper was not only
its contribution to the history of mathematics, but also its contribution
to current research.

Hoping to get some insight into how Siegel spotted and deciphered
these gems among Riemann's fragmentary and disordered personal papers,
we decided to look at the 1932 paper. We also wanted to know how much
of the paper is original to Riemann and whether Siegel needed to fill
any gaps. Although Siegel's paper is cited whenever the Riemann-Siegel
formula is discussed, we were unable to locate an English translation.
Despite our limited knowledge of the German language, we have produced
a translation of the paper as it appears in volume I of Siegel\textquoteright s
collected works.
\end{abstract}
\vspace{10mm}
Comments:
\begin{enumerate}
\item The page and equation numbers in our translation are from the paper
as it appears in Siegel's collected works: \foreignlanguage{ngerman}{\emph{Gesammelte
Abhandlungen}}, Vol. 1, Springer-Verlag, Berlin (1966), pp. 275---310.
\item Our translation is a work in progress. We welcome corrections and
comments on the translation and on mathematical errors we may have
introduce.
\end{enumerate}
\begin{center}
\thispagestyle{empty}\newpage{\LARGE{}On Riemann's Nachlass for
Analytic Number Theory}{\LARGE\par}
\par\end{center}

\begin{center}
\vspace{5mm}
{\large{}Carl Ludwig Siegel (1932)}{\large\par}
\par\end{center}

\vspace{40mm}

In a letter to Weierstrass from the year 1859, Riemann mentioned a
new development of the zeta function which, however, he had not yet
simplified enough for him to be able to include in his published paper
on the theory of prime numbers. Now, since this point of Riemann's
letter had been published by H. Weber in his 1876 edition of Riemann's
works, one could surmise that a detailed review of Riemann's Nachlass
located in the Göttingen University library would yet bring important
hidden formulas of analytic number theory to light.

In fact, Herr Distel, a librarian, already several decades ago discovered
the representation in question of the zeta function in Riemann's papers.
It concerns an asymptotic development which gives the behavior of
the function $\zeta\left(s\right)$ on the critical line $\sigma=\frac{1}{2}$
and, more generally, in each strip $\sigma_{1}\leq\sigma\leq\sigma_{2}$,
expressed for infinitely increasing magnitude of $s$. The principal
term of this development was meanwhile rediscovered, independently
of Riemann, by Hardy and Littlewood in 1920 as a consequence of their
``approximate functional equation''; they use the same methods of
proof as Riemann, namely approximate calculation of an integral by
the saddle-point method. In Riemann, however, there is also a process
for obtaining additional terms of the asymptotic series, which is
based on the beautiful properties of the integral
\[
\varPhi\left(\tau,u\right)=\int\frac{e^{\pi i\tau x^{2}+2\pi iux}}{e^{2\pi ix}-1}dx,
\]
which, incidentally, has also led Kronecker and recently Mordell to
a most elegant derivation of the reciprocity formula for the Gauss
sums. In 1926 Bessel-Hagen noted in a new review of the Riemann papers
another previously unknown
\begin{center}
\thispagestyle{empty}\newpage{}
\par\end{center}

\noindent representation of the zeta function in terms of definite
integrals; in this Riemann was also guided by the properties of $\varPhi\left(\tau,u\right)$.

These two representations of $\zeta\left(s\right)$ may be included
among the most important results in Riemann's number-theoretic Nachlass,
insofar as one does not find them in his published paper. Approaches
to a proof of the so-called ``Riemann hypothesis'' or even to a
proof of the existence of infinitely many zeros of the zeta function
on the critical line are not included in Riemann's papers. On the
conjecture that in the interval $0<t<T$ there lie asymptotically
$\frac{T}{2\pi}\log\frac{T}{2\pi}-\frac{T}{2\pi}$ real zeros of $\zeta\left(\frac{1}{2}+ti\right)$,
Riemann has probably been guided by a heuristic consideration of the
asymptotic series; but even today it is still not clear how one could
prove or disprove this claim. With the help of the asymptotic series
Riemann also calculated some real zeros of $\zeta\left(\frac{1}{2}+ti\right)$
to a better approximation.

In Riemann's notes on the theory of the zeta function there is no
directly publishable material; sometimes disconnected formulas are
found on the same page; often only one side of an equation is written
down; residual estimates and convergence studies are absent, even
at key points. These reasons made a free adaptation of the Riemann
fragments necessary, as has been done in the following.

The legend that Riemann found the results of his mathematical work
through ``highly general'' ideas without requiring the formal tools
of analysis, is probably no longer as widespread as in Klein's lifetime.
Just how strong was Riemann's analytic technique is made particularly
clear by his derivation and manipulation of the asymptotic series
for $\zeta\left(s\right)$.

\section{Evaluation of a Definite Integral.}

Let $u$ be a complex variable. One constructs the integral 
\begin{equation}
\varPhi\left(u\right)=\intop_{\mathclap{0\nwarrow1}}\frac{e^{-\pi ix^{2}+2\pi iux}}{e^{\pi ix}-e^{-\pi ix}}dx,
\end{equation}
extending from $\infty$ to $\infty$ from lower right to upper left
along a line parallel to the bisecting line of the fourth and the
second quadrant, which crosses the real axis between the points 0
and 1. In formula (1) this path of integration is indicated with the
symbol $0\nwarrow1$ placed below the integral sign. \newpage{}

The function $\varPhi\left(u\right)$ is entire. After Riemann it
can be expressed in terms of the exponential function in a simple
way. To prove this, one uses Cauchy's Theorem to obtain two difference
equations for $\varPhi\left(u\right)$:

On the one hand, we have
\begin{align*}
\varPhi\left(u+1\right)-\varPhi\left(u\right) & =\intop_{\mathclap{0\nwarrow1}}e^{-\pi ix^{2}}\frac{e^{2\pi i\left(u+1\right)x}-e^{2\pi iux}}{e^{\pi ix}-e^{-\pi ix}}dx=\intop_{\mathclap{0\nwarrow1}}e^{-\pi ix^{2}+2\pi i\left(u+\frac{1}{2}\right)x}dx\\
 & =e^{\pi i\left(u+\frac{1}{2}\right)^{2}}\intop_{\mathclap{0\nwarrow1}}e^{-\pi i\left(x-u-\frac{1}{2}\right)^{2}}dx=e^{\pi i\left(u+\frac{1}{2}\right)^{2}}\intop_{\mathclap{0\nwarrow1}}e^{-\pi ix^{2}}dx,
\end{align*}
so
\begin{equation}
\varPhi\left(u\right)=\varPhi\left(u+1\right)-e^{\pi i\left(u+\frac{1}{2}\right)^{2}}\intop_{\mathclap{0\nwarrow1}}e^{-\pi ix^{2}}dx.
\end{equation}
On the other hand, when the integration path is indicated by the notation
$-1\nwarrow0$, which arises from the previously used one by a parallel
shift by the vector $-1$, we have
\begin{equation}
1=\intop_{\mathclap{0\nwarrow1}}\frac{e^{-\pi ix^{2}+2\pi iux}}{e^{\pi ix}-e^{-\pi ix}}dx-\intop_{\mathclap{-1\nwarrow0}}\frac{e^{-\pi ix^{2}+2\pi iux}}{e^{\pi ix}-e^{-\pi ix}}dx,
\end{equation}
since the integrand has the residue $\frac{1}{2\pi i}$ at the pole
at $x=0$; so because
\[
\intop_{\mathclap{-1\nwarrow0}}\frac{e^{-\pi ix^{2}+2\pi iux}}{e^{\pi ix}-e^{-\pi ix}}dx=\intop_{\mathclap{0\nwarrow1}}\frac{e^{-\pi i\left(x-1\right)^{2}+2\pi iu\left(x-1\right)}}{e^{\pi i\left(x-1\right)}-e^{-\pi i\left(x-1\right)}}dx=e^{-2\pi iu}\intop_{\mathclap{0\nwarrow1}}\frac{e^{-\pi ix^{2}+2\pi i\left(u+1\right)x}}{e^{\pi ix}-e^{-\pi ix}}dx
\]
equation (3) yields the formula
\begin{equation}
\varPhi\left(u\right)=e^{-2\pi iu}\varPhi\left(u+1\right)+1\,.
\end{equation}

From (2) and (4) one first obtains for $u=0$ the well-known equation
\[
\intop_{\mathclap{0\nwarrow1}}e^{-\pi ix^{2}}dx=e^{\frac{3\pi i}{4}}
\]
and then by elimination of $\varPhi\left(u+1\right)$ the desired
result
\begin{equation}
\intop_{\mathclap{0\nwarrow1}}\frac{e^{-\pi ix^{2}+2\pi iux}}{e^{\pi ix}-e^{-\pi ix}}dx=\frac{1}{1-e^{-2\pi iu}}-\frac{e^{\pi iu^{2}}}{e^{\pi iu}-e^{-\pi iu}}.
\end{equation}
\newpage\noindent Differentiating $n$ times with respect to $u$,
one gets the general formula
\begin{flalign}
 &  &  &  & \intop_{\mathclap{0\nwarrow1}}\frac{e^{-\pi ix^{2}+2\pi iux}}{e^{\pi ix}-e^{-\pi ix}}x^{n}dx & =\left(2\pi i\right)^{-n}D^{n}\frac{e^{\pi iu}-e^{\pi iu^{2}}}{e^{\pi iu}-e^{-\pi iu}} & \left(n=0,1,2,\ldots\right).
\end{flalign}

For the following it is convenient to put (5) in a different form.
One writes $2u+\frac{1}{2}$ instead of $u$ and multiplies (5) by
$e^{-2\pi i\left(u+\frac{1}{2}\right)^{2}+\frac{\pi i}{8}}$; this
gives the equation found by Riemann
\[
\intop_{\mathclap{0\nwarrow1}}\frac{e^{\pi i\left\{ x^{2}-2\left(u+\frac{1}{2}-x\right)^{2}+\frac{1}{8}\right\} }}{e^{2\pi ix}-1}dx=\frac{\cos\left(2\pi u^{2}+\frac{3\pi}{8}\right)}{\cos2\pi u},
\]
which will continue to play an important role.

The integral $\varPhi\left(u\right)$ is a special case of the integral
\begin{equation}
\varPhi\left(\tau,u\right)=\intop_{\mathclap{0\nwarrow1}}\frac{e^{\pi i\tau x^{2}+2\pi iux}}{e^{\pi ix}-e^{-\pi ix}}dx,
\end{equation}
for which two difference equations also suffice. It has been further
studied by Mordell. For each negative rational value of $\tau$ there
is an analogous formula to (5); and from this, through specialization
of $u$, one obtains the reciprocity law for Gauss sums. In his lectures,
Riemann based the transformation theory of the theta functions on
the properties of $\varPhi\left(\tau,u\right)$.

\section{The Asymptotic Formula for the Zeta Function.}

If the real part $\sigma$ of the complex variable $s=\sigma+ti$
is greater than 1 and $m$ is a natural number, then
\[
\zeta\left(s\right)=\sum_{n=1}^{m}n^{-s}+\frac{1}{\varGamma\left(s\right)}\intop_{0}^{\infty}\frac{x^{s-1}e^{-mx}}{e^{x}-1}dx,
\]
or when $C_{1}$ is to be traversed in a loop around the negative
imaginary axis in the positive sense,
\begin{equation}
\zeta\left(s\right)=\sum_{n=1}^{m}n^{-s}+\frac{\left(2\pi\right)^{s}e^{\frac{\pi is}{2}}}{\varGamma\left(s\right)\left(e^{2\pi is}-1\right)}\intop_{C_{1}}\frac{x^{s-1}e^{-2\pi imx}}{e^{2\pi ix}-1}dx.
\end{equation}
This formula is even valid for arbitrary values of $\sigma$. For
now let $\sigma$ be restricted to a fixed interval $\sigma_{1}\leq\sigma\leq\sigma_{2}$,
and let $t>0$. In order to calculate the integral occurring in (8)
asymptotically by the saddle point method for $t\rightarrow\infty$,
one must take the path of integration through the zero\newpage\noindent
point of $D\log\left(x^{s-1}e^{-2\pi imx}\right)$. For this zero
point one obtains from the equation
\[
\frac{s-1}{x}-2\pi im=0
\]
the value
\begin{equation}
\xi=\frac{s-1}{2\pi im}=\frac{t}{2\pi m}+\frac{1-\sigma}{2\pi m}i\,.
\end{equation}
Within the disk centered at $\xi$ of radius $\left|\xi\right|$ the
following development is valid
\begin{align*}
x^{s-1}e^{-2\pi imx} & =\xi^{s-1}e^{-2\pi im\xi}e^{\left(s-1\right)\left\{ -\frac{1}{2}\left(\frac{x-\xi}{\xi}\right)^{2}+\frac{1}{3}\left(\frac{x-\xi}{\xi}\right)^{3}-\cdots\right\} }\\
 & =\xi^{s-1}e^{-2\pi im\xi}e^{-\frac{s-1}{2\xi^{2}}\left(x-\xi\right)^{2}}\left\{ c_{0}+c_{1}\left(x-\xi\right)+c_{2}\left(x-\xi\right)^{2}+\cdots\right\} ;
\end{align*}
and one has in the series
\[
\xi^{s-1}e^{-2\pi im\xi}\sum_{n=0}^{\infty}c_{n}\intop_{C_{1}}\frac{e^{-\frac{s-1}{2\xi^{2}}\left(x-\xi\right)^{2}}}{e^{2\pi ix}-1}\left(x-\xi\right)^{n}dx
\]
a possible asymptotic development of the integral in (8). The integrals
occurring in the series can now all be evaluated by formula (6) of
§ 1, if $\frac{s-1}{2\xi^{2}}$ takes the special value $\pi i$.
For fixed $s$, this is a condition on $m$, which can, in general,
only be satisfied approximately, since $m$ is an integer. This is
why Riemann replaced the saddle point $\xi$ by the neighboring value
$\eta$, which, from the equation
\[
\frac{ti}{2\eta^{2}}=\pi i
\]
is
\begin{equation}
\eta=+\sqrt{\frac{t}{2\pi}}
\end{equation}
and then determined $m$ as suggested by (9) as the greatest integer
less than $\frac{t}{2\pi\eta}$, so
\begin{equation}
m=\left[\eta\right].
\end{equation}

One now introduces the abbreviations
\begin{align}
\tau & =+\sqrt{t}=\eta\sqrt{2\pi}\\
\nonumber \\
\varepsilon & =e^{-\frac{\pi i}{4}}=\dfrac{1-i}{\sqrt{2}}\nonumber \\
\nonumber \\
g\left(x\right) & =x^{s-1}\dfrac{e^{-2\pi imx}}{e^{2\pi ix}-1}\nonumber 
\end{align}
\newpage{}

For now let \foreignlanguage{ngerman}{$\eta$} not be an integer.
The integration path $C_{1}$ will be replaced by the polyline $C_{2}$,
consisting of the two half-lines emanating from the point $\eta-\frac{\varepsilon}{2}\eta$
and containing the points $\eta$ and $-\left(m+\frac{1}{2}\right)$
respectively. With regard to the poles at $\pm1,\pm2,\ldots,\pm m$,
the residue theorem yields 
\[
\intop_{C_{1}}g\left(x\right)dx=\left(e^{\pi is}-1\right)\sum_{n=1}^{m}n^{s-1}+\intop_{C_{2}}g\left(x\right)dx
\]
\begin{equation}
\zeta\left(s\right)=\sum_{n=1}^{m}n^{-s}+\frac{\left(2\pi\right)^{s}}{2\varGamma\left(s\right)\cos\frac{\pi s}{2}}\sum_{n=1}^{m}n^{s-1}+\frac{\left(2\pi\right)^{s}e^{\frac{\pi is}{2}}}{\varGamma\left(s\right)\left(e^{2\pi is}-1\right)}\intop_{C_{2}}g\left(x\right)dx.
\end{equation}
On the left-situated of the two rectilinear elements of $C_{2}$, which
we name $C_{3}$, is now
\[
\arg\,x\geq\arctan\frac{1}{2\sqrt{2}-1}>\left(2\sqrt{2}-1\right)^{-1}-\frac{1}{3}\left(2\sqrt{2}-1\right)^{-3}>\frac{1}{2\sqrt{2}}+\frac{1}{8}
\]
\[
\mathfrak{J}\left(x\right)\leq\frac{\eta}{2\sqrt{2}}
\]
and therefore according to (10) and (11)
\[
\left|x^{s-1}e^{-2\pi imx}\right|\leq\left|x\right|^{\sigma-1}e^{-t\left(\frac{1}{2\sqrt{2}}+\frac{1}{8}\right)+\pi m\frac{\eta}{\sqrt{2}}}\leq\left|x\right|^{\sigma-1}e^{-\frac{t}{8}}
\]
\begin{equation}
\intop_{C_{3}}g\left(x\right)dx=O\left(e^{-\frac{t}{9}}\right),
\end{equation}
uniformly in $\sigma$ for $\sigma_{1}\leq\sigma\leq\sigma_{2}$.

On the right half-line of $C_{2}$ one sets
\begin{flalign*}
 &  & \qquad\qquad x & =\eta+\varepsilon y & \left(y\geq-\frac{\eta}{2}\right);
\end{flalign*}
then we have
\[
\left|x^{s-1}e^{-2\pi imx}\right|=\left|x\right|^{\sigma-1}e^{t\arctan\frac{y}{y+\eta\sqrt{2}}-\pi\sqrt{2}my}.
\]
This is even valid if $y\geq+\frac{\eta}{2}$, so one has for sufficiently
large $t$
\[
t\arctan\frac{y}{y+\eta\sqrt{2}}-\pi\sqrt{2}my\leq\frac{ty}{y+\eta\sqrt{2}}-\pi\sqrt{2}my<ty\left(\frac{1}{y+\eta\sqrt{2}}-\frac{\eta-1}{\sqrt{2}\eta^{2}}\right)
\]
\[
=\frac{ty}{\eta\sqrt{2}}\left(\frac{1}{\eta}-\frac{y}{y+\eta\sqrt{2}}\right)\leq\frac{t}{2\sqrt{2}}\left(\frac{1}{\eta}-\frac{1}{1+2\sqrt{2}}\right)<-\frac{t}{11};
\]
so we have the estimate
\begin{equation}
-\intop_{\frac{\eta}{2}}^{\infty}g\left(x\right)\varepsilon\,dy=O\left(e^{-\frac{t}{11}}\right),
\end{equation}
\newpage\noindent and this in particular is again uniform in $\sigma$
for $\sigma_{1}\leq\sigma\leq\sigma_{2}$. From (14) and (15) one
gets 
\begin{equation}
\intop_{C_{2}}g\left(x\right)dx=\intop_{\eta+\varepsilon\frac{\eta}{2}}^{\eta-\varepsilon\frac{\eta}{2}}g\left(x\right)dx+O\left(e^{-\frac{t}{11}}\right).
\end{equation}

For the asymptotic expansion of the integral on the right side of
(16) one starts from the identity
\begin{equation}
g\left(x\right)=\eta^{s-1}e^{-2\pi im\eta}\frac{e^{-\pi i\left(x-\eta\right)^{2}+2\pi i\left(\eta-m\right)\left(x-\eta\right)}}{e^{2\pi ix}-1}e^{\left(s-1\right)\log\left(1+\frac{x-\eta}{\eta}\right)-2\pi i\eta\left(x-\eta\right)+\pi i\left(x-\eta\right)^{2}}
\end{equation}
For $\left|x-\eta\right|<\eta$, the last of the right-side factors
can be developed in a series in powers of $x-\eta$, whose coefficients
are to be investigated in more detail. With the definition of $\tau$
stated in (12), one sets
\begin{flalign}
 &  &  &  & e^{\left(s-1\right)\log\left(1+\frac{z}{\tau}\right)-i\tau z+\frac{i}{2}z^{2}} & =\sum_{n=0}^{\infty}a_{n}z^{n}=w\left(z\right) & \left(\left|z\right|<\tau\right).
\end{flalign}
From the differential equation
\[
\left(z+\tau\right)\frac{dw}{dz}+\left(1-\sigma-iz^{2}\right)w=0
\]
there results the recursion formula 
\begin{flalign}
 &  &  &  & \left(n+1\right)\tau a_{n+1} & =-\left(n+1-\sigma\right)a_{n}+ia_{n-2} & \left(n=2,3,\ldots\right),
\end{flalign}
which is also correct for $n=0,1$, when one sets $a_{-2}=0,a_{-1}=0$.
If one adds the equation $a_{0}=1$ to this, then $a_{1},a_{2}\ldots$
are determined by virtue of (19); in particular, $a_{n}$ is a polynomial
of the $n^{\mathrm{th}}$ degree in $\tau^{-1}$, that does not contain
the powers $\tau^{-k}$ for $k=0,1,\ldots n-2\left[\frac{n}{3}\right]-1$.
Consequently 
\[
a_{n}=O\left(t^{-\frac{n}{2}+\left[\frac{n}{3}\right]}\right)
\]
is true uniformly for $\sigma_{1}\leq\sigma\leq\sigma_{2}$, but somewhat
non-uniformly in $n$.

To estimate the tail of the power series $w\left(z\right)$, one uses
the representation
\begin{equation}
r_{n}\left(z\right)=\sum_{k=n}^{\infty}a_{k}z^{k}=\frac{1}{2\pi i}\intop_{C}\frac{w\left(u\right)z^{n}}{u^{n}\left(u-z\right)}du;
\end{equation}
where $C$ is a curve lying in the convergence circle, which loops
positively around the points 0 and $z$ once each. By (18) we have
\begin{align*}
\log w\left(u\right) & =\left(\sigma-1+i\tau^{2}\right)\log\left(1+\frac{u}{\tau}\right)-i\tau u+\frac{i}{2}u^{2}\\
\\
 & =\left(\sigma-1\right)\log\left(1+\frac{u}{\tau}\right)+iu^{2}\sum_{k=1}^{\infty}\frac{\left(-1\right)^{k-1}}{k+2}\left(\frac{u}{\tau}\right)^{k};
\end{align*}
\newpage\noindent thus, valid in the circle $\left|u\right|\leq\frac{3}{5}\tau$,
we have the estimate
\begin{equation}
\mathfrak{R}\log w\left(u\right)\leq\left|\sigma-1\right|\log\frac{5}{2}+\frac{5}{6}\frac{\left|u\right|}{\tau}\left|u\right|^{2}.
\end{equation}
In (20) set $\left|z\right|\leq\frac{4}{7}\tau$ and let $C$ be a
circle around $u=0$ with a radius $\varrho_{n}$, which, for now,
will only be subject to the condition
\begin{equation}
\frac{21}{20}\left|z\right|\leq\varrho_{n}\leq\frac{3}{5}\tau
\end{equation}
From (20), (21), (22) then follows uniformly in $\sigma$ and $n$
the estimate
\begin{equation}
r_{n}\left(z\right)=O\left(\left|z\right|^{n}\varrho_{n}^{-n}e^{\frac{5}{6\tau}\varrho_{n}^{3}}\right).
\end{equation}
The function $\varrho^{-n}e^{\frac{5}{6\tau}\varrho^{3}}$ of $\varrho$
has its minimum $\left(\frac{5e}{2n\tau}\right)^{\frac{n}{3}}$ for
$\varrho=\left(\frac{2n\tau}{5}\right)^{\frac{1}{3}}$. By (22) the
choice $\varrho_{n}=\varrho$ is permitted, in the case that 
\[
\frac{21}{20}\left|z\right|\leq\left(\frac{2n\tau}{5}\right)^{\frac{1}{3}}\leq\frac{3}{5}\tau
\]
Consequently, we have 
\begin{flalign}
 &  &  &  &  &  &  & r_{n}\left(z\right)=O\left(\left|z\right|^{n}\left(\frac{5e}{2n\tau}\right)^{\frac{n}{3}}\right) & \left(n\leq\frac{27}{50}t,\left|z\right|\leq\frac{20}{21}\left(\frac{2n\tau}{5}\right)^{\frac{1}{3}}\right).
\end{flalign}
For $\left|z\right|\leq\frac{4}{7}\tau$, by (22), the choice $\varrho_{n}=\frac{21}{20}\left|z\right|$
is also admissible; then (23) yields the relation 
\begin{flalign}
 &  &  &  & r_{n}\left(z\right) & =O\left(\left(\frac{20}{21}\right)^{n}e^{\frac{5}{6\tau}\left(\frac{21}{20}\left|z\right|\right)^{3}}\right)=O\left(e^{\frac{14}{29}\left|z\right|^{2}}\right) & \left(\left|z\right|\leq\frac{\tau}{2}\right).
\end{flalign}

By (17) and (18) it follows that
\begin{equation}
\intop_{\eta+\varepsilon\frac{\eta}{2}}^{\eta-\varepsilon\frac{\eta}{2}}g\left(x\right)dx=\eta^{s-1}e^{-2\pi im\eta}\intop_{\eta+\varepsilon\frac{\eta}{2}}^{\eta-\varepsilon\frac{\eta}{2}}\frac{e^{-\pi i\left(x-\eta\right)^{2}+2\pi i\left(\eta-m\right)\left(x-\eta\right)}}{e^{2\pi ix}-1}w\left(\sqrt{2\pi}\left(x-\eta\right)\right)dx\,.
\end{equation}
To determine the error that results if one replaces in this equation
$w\left(\sqrt{2\pi}\left(x-\eta\right)\right)$ by the partial sum
$\sum_{k=0}^{n-1}a_{k}\left(2\pi\right)^{\frac{k}{2}}\left(x-\eta\right)^{k}$,
one must examine the integral
\begin{equation}
J_{n}=\intop_{\eta+\varepsilon\frac{\eta}{2}}^{\eta-\varepsilon\frac{\eta}{2}}\frac{e^{-\pi i\left(x-\eta\right)^{2}+2\pi i\left(\eta-m\right)\left(x-\eta\right)}}{e^{2\pi ix}-1}r_{n}\left(\sqrt{2\pi}\left(x-\eta\right)\right)dx
\end{equation}
\newpage\noindent From now on let $n\leq\frac{5}{16}t$. The neighborhoods
of the poles $x=m,m+1$ of the integrand will thereby be avoided,
if one replaces the part of the path of integration that lies within
the disks $\left|x-m\right|\leq\frac{1}{2\sqrt{\pi}}$ or $\left|x-m-1\right|\leq\frac{1}{2\sqrt{\pi}}$
respectively by the corresponding arc. As in (24), the integration
over the circular arc provides only the contribution $O\left(\left(\frac{5e}{2n\tau}\right)^{\frac{n}{3}}\right)$
to $J_{n}$. On the remainder of the path of integration one has $-\pi i\left(x-\eta\right)^{2}=-\pi\left|x-\eta\right|^{2}$.
One sets 
\[
\frac{20}{21}\left(\frac{2n\tau}{5}\right)^{\frac{1}{3}}=\lambda
\]
and considers (24) for $\left|x-\eta\right|\leq\frac{1}{\sqrt{2\pi}}$,
or on the other hand (25) for $\frac{1}{\sqrt{2\pi}}\leq\left|x-\eta\right|\leq\frac{\eta}{2}$;
then there follows
\begin{align*}
J_{n} & =O\left\{ \left(\frac{5e}{2n\tau}\right)^{\frac{n}{3}}\intop_{0}^{\lambda}e^{-\frac{1}{2}v^{2}+\sqrt{2\pi}v}v^{n}dv+\intop_{\lambda}^{\frac{\tau}{2}}e^{-\frac{1}{58}v^{2}+\sqrt{2\pi}v}dv\right\} \\
\\
 & =O\left\{ \left(\frac{5e}{2n\tau}\right)^{\frac{n}{3}}e^{\sqrt{2\pi n}}2^{\frac{n}{2}}\varGamma\left(\frac{n+1}{2}\right)+e^{-\frac{1}{59}\lambda^{2}}\right\} =O\left\{ \left(\frac{25n}{4et}\right)^{\frac{n}{6}}e^{\sqrt{2\pi n}}+e^{-\frac{1}{59}\lambda^{2}}\right\} .
\end{align*}
A simple calculation shows that for $n\leq2\cdot10^{-8}\,t$ the second
$O$-term is exceeded by the first. This yields the estimate 
\begin{flalign}
 &  &  &  &  & J_{n}=O\left(\left(\frac{3n}{t}\right)^{\frac{n}{6}}\right) & \left(n\leq2\cdot10^{-8}\,t\right)
\end{flalign}
uniformly in $\sigma$ and $n$.

From (16), (18), (26), (27), (28) now follows
\[
\intop_{C_{2}}g\left(x\right)dx
\]
\[
=\eta^{s-1}e^{-2\pi im\eta}\left\{ \sum_{k=0}^{n-1}a_{k}\left(2\pi\right)^{\frac{k}{2}}\intop_{\eta+\varepsilon\frac{\eta}{2}}^{\eta-\varepsilon\frac{\eta}{2}}\frac{e^{-\pi i\left(x-\eta\right)^{2}+2\pi i\left(\eta-m\right)\left(x-\eta\right)}}{e^{2\pi ix}-1}\left(x-\eta\right)^{k}dx+O\left(\left(\frac{3n}{t}\right)^{\frac{n}{6}}\right)\right\} .
\]
Integrating on the right side instead of $\eta+\varepsilon\frac{\eta}{2}$
to $\eta-\varepsilon\frac{\eta}{2}$ over the full line of $\eta+\varepsilon\,\infty$
to $\eta-\varepsilon\,\infty$, then, since $n\leq2\cdot10^{-8}\,t$
{[}Red Herring?{]}, the value of the integral $\left[\textrm{Translater's Note: }\intop_{C_{2}}g\left(x\right)dx\right]$
only changes by $O\left(e^{-\frac{t}{8}+\pi\eta}\left(\frac{\eta}{2}\right)^{k}\right)$;
on the other hand, according to (24)\newpage
\begin{flalign*}
 &  &  &  &  &  &  & a_{k}=\left(r_{k}-r_{k+1}\right)z^{-k}=O\left(\left(\frac{5e}{2k\tau}\right)^{\frac{k}{3}}\right) & \left(k=1,2,\ldots,n-1\right),
\end{flalign*}
 so
\[
\sum_{k=0}^{n-1}\left|a_{k}\right|e^{-\frac{t}{8}+\pi\eta}\left(\frac{\tau}{2}\right)^{k}=O\left(e^{-\frac{t}{8}+\pi\eta}\left(\frac{5et}{16n}\right)^{\frac{n}{3}}\right)=O\left(\left(\frac{3n}{t}\right)^{\frac{n}{6}}\right).
\]
Finally, one replaces the integration variable $x$ by $x+m$, so
that
\begin{equation}
\intop_{C_{2}}g\left(x\right)dx
\end{equation}
\[
=\left(-1\right)^{m}e^{-\frac{\pi i}{8}}\eta^{s-1}e^{-\pi i\eta^{2}}\left\{ \sum_{k=0}^{n-1}a_{k}\left(2\pi\right)^{\frac{k}{2}}\intop_{\mathclap{0\nwarrow1}}\frac{e^{\pi i\left(x^{2}-2\left(x+m-\eta\right)^{2}+\frac{1}{8}\right)}}{e^{2\pi ix}-1}\left(x+m-\eta\right)^{k}dx+O\left(\left(\frac{3n}{t}\right)^{\frac{n}{6}}\right)\right\} .
\]

According to the results of § 1 the integral 
\begin{equation}
\intop_{\mathclap{0\nwarrow1}}\frac{e^{\pi i\left(x^{2}-2\left(x-\frac{u}{\sqrt{2\pi}}-\frac{1}{2}\right)^{2}+\frac{1}{8}\right)}}{e^{2\pi ix}-1}dx=F\left(u\right)
\end{equation}
has the value
\[
F\left(u\right)=\frac{\cos\left(u^{2}+\frac{3\pi}{8}\right)}{\cos\left(\sqrt{2\pi}u\right)}.
\]
To also express the right integral occurring in (29) in an elementary
way for $k>0$, Riemann forms from (30) the equation 
\[
F\left(\delta+u\right)e^{iu^{2}}=\intop_{\mathclap{0\nwarrow1}}\frac{e^{\pi i\left(x^{2}-2\left(x-\frac{\delta}{\sqrt{2\pi}}-\frac{1}{2}\right)^{2}+\frac{1}{8}\right)}}{e^{2\pi ix}-1}e^{2\sqrt{2\pi}i\left(x-\frac{\delta}{\sqrt{2\pi}}-\frac{1}{2}\right)u}dx,
\]
from which he obtains, by expansion in powers of $u$, the formula
\begin{equation}
\intop_{\mathclap{0\nwarrow1}}\frac{e^{\pi i\left(x^{2}-2\left(x-\frac{\delta}{\sqrt{2\pi}}-\frac{1}{2}\right)^{2}+\frac{1}{8}\right)}}{e^{2\pi ix}-1}\left(x-\frac{\delta}{\sqrt{2\pi}}-\frac{1}{2}\right)^{k}dx
\end{equation}
\begin{flalign*}
 &  &  &  &  & =2^{-k}\left(2\pi\right)^{-\frac{k}{2}}k!\sum_{r=0}^{\left[\frac{k}{2}\right]}\frac{i^{r-k}}{r!\left(k-2r\right)!}F^{\left(k-2r\right)}\left(\delta\right) & \left(k=0,1,2,\ldots\right)
\end{flalign*}

Now, from (13), (29), (31) there follows the development
\begin{align}
\zeta\left(s\right) & =\sum_{l=1}^{m}l^{-s}+\frac{\left(2\pi\right)^{s}}{2\varGamma\left(s\right)\cos\frac{\pi s}{2}}\sum_{l=1}^{m}l^{s-1}\\
\nonumber \\
 & +\left(-1\right)^{m-1}\frac{\left(2\pi\right)^{\frac{s+1}{2}}}{\varGamma\left(s\right)}t^{\frac{s-1}{2}}e^{\frac{\pi is}{2}-\frac{ti}{2}-\frac{\pi i}{8}}S\nonumber 
\end{align}
\newpage\noindent with
\begin{equation}
S=\sum_{0\leq2r\leq k\leq n-1}\frac{2^{-k}i^{r-k}k!}{r!\left(k-2r\right)!}\,a_{k}\,F^{\left(k-2r\right)}\left(\delta\right)+O\left(\left(\frac{3n}{t}\right)^{\frac{n}{6}}\right)
\end{equation}
\[
n\leq2\cdot10^{-8}\,t,\quad\delta=\sqrt{t}-\left(m+\frac{1}{2}\right)\sqrt{2\pi},\quad F\left(u\right)=\frac{\cos\left(u^{2}+\frac{3\pi}{8}\right)}{\cos\left(\sqrt{2\pi}u\right)};
\]
and where the coefficients $a_{k}$ are determined by the recurrence
formula (19). This development is asymptotic, and in particular is
uniform for $\sigma_{1}\leq\sigma\leq\sigma_{2}$, because the remainder
term in (33) is indeed, for every fixed $n$, of the order $t^{-\frac{n}{6}}$
uniformly in $\sigma$. From the known asymptotic series of analysis,
(32) itself differs by the appearance of the integer $m$, which causes
the individual terms of the expansion to depend on $t$ discontinuously.
The assumption made in the proof that $\sqrt{\frac{t}{2\pi}}$ is
not an integer, can easily be subsequently eliminated, because one
can indeed make in (32) the right-hand boundary crossing at any integer
values of $\sqrt{\frac{t}{2\pi}}$.

If we choose in (33) the special value $n=2\cdot10^{-8}\,t$, the
error term is $O\left(10^{-10^{-8}\,t}\right)$, which therefore goes
exponentially to 0 with increasing $t$. For practical purposes, this
estimate of the error term is not useful because of the small exponent
factor $10^{-8}$; finer estimates show that $10^{-8}$ can be replaced
by a considerably larger number. It would be of interest to find the
exact order of growth of the error as a function of $n$; but this
is not trivial, since it does not converge to 0 as $n$ increases
for fixed $t$.

Because of the particular importance of the case $\sigma=\frac{1}{2}$,
it is expedient to multiply (32) by the function $e^{\vartheta i}$
defined by
\begin{equation}
e^{\vartheta i}=\pi^{\frac{1}{4}-\frac{s}{2}}\sqrt{\frac{\varGamma\left(\frac{s}{2}\right)}{\varGamma\left(\frac{1-s}{2}\right)}}
\end{equation}
where we mean by $\vartheta$ those values in the uniquely determined
branch in the plane cut from 0 to $-\infty$ and from 1 to $+\infty$,
which vanishes for $s=\frac{1}{2}$. Then, on the critical line $\sigma=\frac{1}{2}$,
$\vartheta=\mathrm{\arg\left(\pi^{-\frac{s}{2}}\varGamma\left(\frac{s}{2}\right)\right)}$
and $e^{\vartheta i}\zeta\left(s\right)$ is real. Following (32),
for $\sigma_{1}\leq\sigma\leq\sigma_{2}$ \newpage
\begin{align}
e^{\vartheta i}\zeta\left(s\right) & =2\sum_{l=1}^{m}\frac{\cos\left(\vartheta+i\left(s-\frac{1}{2}\right)\log l\right)}{\sqrt{l}}\\
\nonumber \\
 & +\left(-1\right)^{m-1}\left(\frac{t}{2\pi}\right)^{\frac{\sigma-1}{2}}e^{\left(\frac{t}{2}\log\frac{t}{2\pi}-\frac{t}{2}-\frac{\pi}{8}-\vartheta\right)i}S,\nonumber 
\end{align}
 where $S$ is defined by (33). Each $a_{k}$ included in the finite
sum in $S$ is a polynomial in $\tau^{-1}$. It therefore follows
from (33), for each fixed $n$ and $t\rightarrow\infty$, that by
re-ordering in powers of $\tau^{-1}$ we have the relationship 
\[
S=\sum_{k=0}^{n-1}A_{k}\tau^{-k}+O\left(\tau^{-n}\right),
\]
where the coefficients $A_{0},A_{1},\ldots,A_{n-1}$ are homogeneous
linear combinations of a finite number of the derivatives $F\left(\delta\right),F'\left(\delta\right),\ldots$
. The explicit calculation of the $A_{k}$ with the help of (33) and
the recurrence formula for $a_{k}$ is rather cumbersome; Riemann
simplified this by the following artifice. Substituting 
\[
F\left(\delta+x\right)e^{ix^{2}}=\sum_{k=0}^{\infty}b_{k}x^{k},
\]
then
\begin{equation}
S\sim\sum_{k=0}^{\infty}\left(2i\right)^{-k}k!a_{k}b_{k}
\end{equation}
is a full asymptotic series, and the required quantity $A_{k}$ is
the coefficient of $\tau^{-k}$, which arises by ordering the series
in powers of $\tau^{-1}$. The sum on the right hand side of (36)
is but the constant term in the series of positive and negative powers
of $x$, that results through formal multiplication of the convergent
power series
\[
F\left(\delta+\frac{1}{x}\right)e^{ix^{-2}}=\sum_{k=0}^{\infty}b_{k}x^{-k}
\]
with the divergent power series
\begin{equation}
y=\sum_{k=0}^{\infty}\left(2i\right)^{-k}k!\,a_{k}\,x^{k}.
\end{equation}
Since the fixed power $\tau^{-k}$ only occurs in a finite number
of coefficients $a_{0},a_{1},a_{2},\ldots$, then the following method
for the calculation of $A_{k}$ is legitimate: Construct by formal
multiplication of $e^{ix^{-2}}$ and $y$ the term
\begin{equation}
z=e^{ix^{-2}}y=\sum_{n=-\infty}^{+\infty}d_{n}x^{n}
\end{equation}
\newpage\noindent and from this, by ordering in powers of $\tau^{-1}$,
the series $\sum_{k=0}^{\infty}B_{k}\tau^{-k}$; then $A_{k}$ is
the constant term in $F\left(\delta+\frac{1}{x}\right)B_{k}$. For
the calculation of this constant term, the negative powers of $x$
occurring in the $B_{k}$ are irrelevant, so one only has to determine
the polynomial part of $B_{k}$.

If, for brevity, one sets
\begin{flalign*}
 &  &  &  & \left(2i\right)^{-k}k!a_{k} & =c_{k} & \left(k=0,1,2,\ldots\right),
\end{flalign*}
then, by (19)
\begin{flalign*}
 &  &  &  & \tau c_{n+1} & =i\frac{n+1-\sigma}{2}c_{n}-\frac{n\left(n-1\right)}{8}c_{n-2} & \left(n=0,1,2,\ldots\right)
\end{flalign*}
with $c_{-2}=0,c_{-1}=0,c_{0}=1$, and hence the power series (37)
formally satisfies the differential equation 
\[
\tau\left(y-1\right)=\frac{i}{2}x^{\sigma+1}D\left(x^{1-\sigma}y\right)-\frac{1}{8}x^{3}D^{2}\left(x^{2}y\right).
\]
It follows as a differential equation for the series (38)
\begin{equation}
\left\{ \tau+\frac{1}{2x}+i\left(\frac{\sigma}{2}-\frac{1}{4}\right)x\right\} z+\frac{1}{8}x^{3}D^{2}\left(x^{2}z\right)=\tau e^{ix^{-2}}.
\end{equation}
If now, ordered in powers of $\tau^{-1}$,
\[
z=\sum_{n=0}^{\infty}B_{n}\tau^{-n},
\]
then it follows from (39)
\[
B_{0}=e^{ix^{-2}}
\]
and also the recursion formula
\begin{flalign*}
 &  &  &  &  &  & B_{n+1} & =\left(i\frac{1-2\sigma}{4}x-\frac{1}{2x}\right)B_{n}-\frac{1}{8}x^{3}D^{2}\left(x^{2}B_{n}\right) & \left(n=0,1,2,\ldots\right).
\end{flalign*}
If one now sets
\begin{flalign*}
 &  &  &  & B_{n} & =\sum_{k=-\infty}^{3n}a_{k}^{\left(n\right)}x^{k} & \left(n=0,1,2,\ldots\right),
\end{flalign*}
then one has
\begin{flalign*}
 &  &  &  & a_{k}^{\left(0\right)} & =0 & \left(k\neq0,-2,-4,-6,\ldots\right)\\
 &  &  &  & a_{-2k}^{\left(0\right)} & =\frac{i^{k}}{k!} & \left(k=0,1,2,3,\ldots\right)
\end{flalign*}
\begin{multline}
\qquad\qquad\qquad a_{k}^{\left(n+1\right)}=i\frac{1-2\sigma}{4}a_{k-1}^{\left(n\right)}-\frac{1}{2}a_{k+1}^{\left(n\right)}-\frac{\left(k-1\right)\left(k-2\right)}{8}a_{k-3}^{\left(n\right)}\\
\left(n=0,1,2,\ldots;k=0,\pm1,\pm2,\ldots\right).\!\!\!\!
\end{multline}
With the help of these recursion formulas for calculating the $a_{k}^{\left(n\right)}$,
then the $A_{n}$ can be explicitly specified, namely
\begin{equation}
A_{n}=\sum_{k=0}^{3n}\frac{a_{k}^{\left(n\right)}}{k!}F^{\left(k\right)}\left(\delta\right);
\end{equation}
\newpage\noindent and this yields
\[
S\sim\sum\frac{a_{k}^{\left(n\right)}}{k!}F^{\left(k\right)}\left(\delta\right)\tau^{-n},
\]
where $n$ takes the values $0,1,2,\ldots$ and $k$ the values of
$0,1,2,\ldots,3n$.

The recursion formula (40) is easiest for $\sigma=\frac{1}{2}$. For
this special case one easily obtains
\begin{align*}
B_{0} & =1+\cdots\\
\\
B_{1} & =-\tfrac{1}{2^{2}}x^{3}+\cdots\\
\\
B_{2} & =\tfrac{5}{2^{3}}x^{6}+\tfrac{1}{2^{3}}x^{2}+\tfrac{i}{2^{4}\cdot3}+\cdots\\
\\
B_{3} & =-\tfrac{5\cdot7}{2^{3}}x^{9}-\tfrac{1}{2}x^{5}-\tfrac{i}{2^{6}\cdot3}x^{3}-\tfrac{1}{2^{4}}x+\cdots\\
\\
B_{4} & =\tfrac{5^{2}\cdot7\cdot11}{2^{5}}x^{12}+\tfrac{7\cdot11}{2^{4}}x^{8}+\tfrac{5i}{2^{7}\cdot3}x^{6}+\tfrac{19}{2^{6}}x^{4}+\tfrac{i}{3\cdot2^{7}}x^{2}+\tfrac{11\cdot13}{2^{9}\cdot3^{2}}+\cdots;
\end{align*}
where the omitted summands contain only negative powers of $x$. Consequently,
for $\sigma=\frac{1}{2}$ 
\begin{equation}
\left\{ \begin{aligned}A_{0} & =F\left(\delta\right)\\
\\
A_{1} & =-\tfrac{1}{2^{3}\cdot3}F^{\left(3\right)}\left(\delta\right)\\
\\
A_{2} & =\tfrac{1}{2^{7}\cdot3^{2}}F^{\left(6\right)}\left(\delta\right)+\tfrac{1}{2^{4}}F^{\left(2\right)}\left(\delta\right)+\tfrac{i}{2^{4}\cdot3}F\left(\delta\right)\\
\\
A_{3} & =-\tfrac{1}{2^{10}\cdot3^{4}}F^{\left(9\right)}\left(\delta\right)-\tfrac{1}{2^{4}\cdot3\cdot5}F^{\left(5\right)}\left(\delta\right)-\tfrac{i}{2^{7}\cdot3^{2}}F^{\left(3\right)}\left(\delta\right)-\tfrac{1}{2^{4}}F^{\left(1\right)}\left(\delta\right)\\
\\
A_{4} & =\tfrac{1}{2^{15}\cdot3^{5}}F^{\left(12\right)}\left(\delta\right)+\tfrac{11}{2^{11}\cdot3^{2}\cdot5}F^{\left(8\right)}\left(\delta\right)+\tfrac{i}{2^{11}\cdot3^{3}}F^{\left(6\right)}\left(\delta\right)\\
\\
 & +\tfrac{19}{2^{9}\cdot3}F^{\left(4\right)}\left(\delta\right)+\tfrac{i}{2^{8}\cdot3}F^{\left(2\right)}\left(\delta\right)+\tfrac{11\cdot13}{2^{9}\cdot3^{2}}F\left(\delta\right);
\end{aligned}
\right.
\end{equation}
and thus $S$ is determined in the case of $\sigma=\frac{1}{2}$ to
an error on the order of $\tau^{-5}$.

The asymptotic development (35) can yet be simplified somewhat, if
one develops the quantity $\vartheta$ in the second term on the right
side asymptotically with the help of Stirling's series. For this purpose,
Riemann considers the formula 
\[
\log\varGamma\left(\frac{1}{4}+\frac{ti}{2}\right)
\]
\begin{flalign*}
 &  &  &  &  & =\left(\frac{ti}{2}-\frac{1}{4}\right)\log\frac{ti}{2}-\frac{ti}{2}+\log\sqrt{2\pi}+\frac{1}{4}\intop_{0}^{\infty}\left(\frac{4e^{3x}}{e^{4x}-1}-\frac{1}{x}-1\right)\frac{e^{-2tix}}{x}dx & \left(t>0\right),
\end{flalign*}
resulting from the known Binet's integral representation of $\log\varGamma\left(s\right)$
by a simple transformation. Because of the identity\newpage
\[
\frac{4e^{3x}}{e^{4x}-1}=\frac{1}{\cosh x}+\frac{1}{\sinh x}
\]
 it follows by separation of real and imaginary parts 
\begin{eqnarray*}
\log\left|\varGamma\left(\frac{1}{4}+\frac{ti}{2}\right)\right| & = & -\frac{\pi}{4}t-\frac{1}{4}\log\frac{t}{2}+\log\sqrt{2\pi}+\frac{1}{4}\intop_{0}^{\infty}\left(\frac{1}{\cos x}-1\right)\frac{e^{-2tx}}{x}dx\\
\\
 & - & \frac{1}{4}\log\left(1+e^{-2\pi t}\right)\\
\\
\arg\varGamma\left(\frac{1}{4}+\frac{ti}{2}\right) & = & \frac{t}{2}\log\frac{t}{2}-\frac{t}{2}-\frac{\pi}{8}+\frac{1}{4}\intop_{0}^{\infty}\left(\frac{1}{\sin x}-\frac{1}{x}\right)\frac{e^{-2tx}}{x}dx\\
\\
 & + & \frac{1}{2}\arctan e^{-\pi t},
\end{eqnarray*}
where the integrals are to be understood as Cauchy principal values
due to the poles at $k\frac{\pi}{2}\left(k=1,2,\ldots\right)$. If
one sets
\begin{flalign*}
 &  &  & \frac{1}{\cos x}=\sum_{n=0}^{\infty}\frac{E_{n}}{\left(2n\right)!}x^{2n} & \left(\left|x\right|<\frac{\pi}{2}\right)\\
 &  &  & \frac{x}{\sin x}=\sum_{n=0}^{\infty}\frac{F_{n}}{\left(2n\right)!}x^{2n} & \left(\left|x\right|<\pi\right),
\end{flalign*}
then $E_{0}=1,E_{1}=1,E_{2}=5,E_{3}=61,F_{0}=1,F_{1}=\frac{1}{3},F_{2}=\frac{7}{15},F_{3}=\frac{31}{21}$
and generally
\begin{flalign*}
 &  &  &  & E_{n} & -\dbinom{2n}{2}E_{n-1}+\dbinom{2n}{4}E_{n-2}-\cdots+\left(-1\right)^{n}E_{0}=0 & \left(n=1,2,3,\ldots\right)\\
\\
 &  &  &  & \dbinom{2n+1}{1}F_{n} & -\dbinom{2n+1}{3}F_{n-1}+\dbinom{2n+1}{5}F_{n-2}-\cdots+\left(-1\right)^{n}F_{0}=0\\
 &  &  &  &  &  & \left(n=1,2,3,\ldots\right).
\end{flalign*}
This provides in a well-known manner the asymptotic series
\begin{equation}
\log\left|\varGamma\left(\frac{1}{4}+\frac{ti}{2}\right)\right|\sim-\frac{\pi}{4}t-\frac{1}{4}\log\frac{t}{2}+\log\sqrt{2\pi}+\frac{1}{8}\sum_{n=1}^{\infty}\frac{E_{n}}{n}\left(2t\right)^{-2n}
\end{equation}
\[
\arg\varGamma\left(\frac{1}{4}+\frac{ti}{2}\right)\sim\frac{t}{2}\log\frac{t}{2}-\frac{t}{2}-\frac{\pi}{8}+\frac{1}{8}\sum_{n=1}^{\infty}\frac{F_{n}}{n\left(2n-1\right)}\left(2t\right)^{1-2n}.
\]

On $\sigma=\frac{1}{2}$ there now holds $\vartheta=-\frac{t}{2}\log\pi+\arg\varGamma\left(\frac{1}{4}+\frac{ti}{2}\right)$
and consequently

\[
\frac{t}{2}\log\frac{t}{2\pi}-\frac{t}{2}-\frac{\pi}{8}-\vartheta\sim-\frac{1}{8}\sum_{n=1}^{\infty}\frac{F_{n}}{n\left(2n-1\right)}\left(2t\right)^{1-2n}
\]
\[
e^{\left(\frac{t}{2}\log\frac{t}{2\pi}-\frac{t}{2}-\frac{\pi}{8}-\vartheta\right)i}=1-\frac{i}{2^{4}\cdot3}t^{-1}-\frac{1}{2^{9}\cdot3^{2}}t^{-2}+O\left(t^{-3}\right).
\]
\newpage\noindent Taking into account (42) we now have, as a definitive
form of asymptotic series for $\zeta\left(s\right)$ for $\sigma=\frac{1}{2}$,
the equation
\begin{equation}
e^{\vartheta i}\zeta\left(\frac{1}{2}+ti\right)=2\sum_{n=1}^{m}\frac{\cos\left(\vartheta-t\log n\right)}{\sqrt{n}}+\left(-1\right)^{m-1}\left(\frac{t}{2\pi}\right)^{-\frac{1}{4}}R
\end{equation}
where
\[
\vartheta=-\frac{t}{2}\log\pi+\arg\varGamma\left(\frac{1}{4}+\frac{ti}{2}\right)
\]
\[
m=\left[\sqrt{\frac{t}{2\pi}}\;\;\right]
\]
\[
R\sim C_{0}+C_{1}t^{-\frac{1}{2}}+C_{2}t^{-1}+C_{3}t^{-\frac{3}{2}}+C_{4}t^{-2}+\cdots
\]
\begin{equation}
\left\{ \begin{aligned}C_{0} & =F\left(\delta\right)\\
\\
C_{1} & =-\tfrac{1}{2^{3}\cdot3}F^{\left(3\right)}\left(\delta\right)\\
\\
C_{2} & =\tfrac{1}{2^{4}}F^{\left(2\right)}\left(\delta\right)+\tfrac{1}{2^{7}\cdot3^{2}}F^{\left(6\right)}\left(\delta\right)\\
\\
C_{3} & =-\tfrac{1}{2^{4}}F^{\left(1\right)}\left(\delta\right)-\tfrac{1}{2^{4}\cdot3\cdot5}F^{\left(5\right)}\left(\delta\right)-\tfrac{1}{2^{10}\cdot3^{4}}F^{\left(9\right)}\left(\delta\right)\\
\\
C_{4} & =\tfrac{1}{2^{5}}F\left(\delta\right)+\tfrac{19}{2^{9}\cdot3}F^{\left(4\right)}\left(\delta\right)+\tfrac{11}{2^{11}\cdot3^{2}\cdot5}F^{\left(8\right)}\left(\delta\right)+\tfrac{1}{2^{15}\cdot3^{5}}F^{\left(12\right)}\left(\delta\right)
\end{aligned}
\right.
\end{equation}
\[
\begin{aligned}F\left(x\right) & =\frac{\cos\left(x^{2}+\frac{3\pi}{8}\right)}{\cos\left(\sqrt{2\pi}x\right)}\\
\delta & =\sqrt{t}-\left(m+\frac{1}{2}\right)\sqrt{2\pi},
\end{aligned}
\]
and one essentially finds the above in Riemann. The only new contribution
in the foregoing is the remainder estimate.

One may now drop the condition $\sigma=\frac{1}{2}$ and the restriction
of $\sigma$ only to the interval $\sigma_{1}\leq\sigma\leq\sigma_{2}$,
and still can use the asymptotic development (44); then one must just
understand by $t$ the complex number $-i\left(s-\frac{1}{2}\right)$
and by $m$ the integer $\left[\sqrt{\frac{\left|t\right|}{2\pi}}\;\right]$,
while the meaning of $\vartheta$ is still defined by (34). The needed
extensions to prove this assertion can be obtained without difficulty
from the foregoing derivation of (44).

The asymptotic series $R$ is a homogeneous linear combination of
the quantities $F\left(\delta\right)$, $F'\left(\delta\right)$,
$F''\left(\delta\right)$,\foreignlanguage{ngerman}{$\,\ldots$};
by rearrangement there arises from it an expression of the form
\[
D_{0}^{*}F\left(\delta\right)+D_{1}^{*}F'\left(\delta\right)+D_{2}^{*}F''\left(\delta\right)+\cdots,
\]
\newpage\noindent where each \foreignlanguage{ngerman}{$D_{n}^{\ast}$}
is a power series in $\tau^{-1}$. These power series are divergent;
it raises the question, whether they are asymptotic developments of
certain analytic functions $D_{0},D_{1},D_{2},\ldots$ and whether
the series 
\begin{equation}
D_{0}F\left(\delta\right)+D_{1}F'\left(\delta\right)+D_{2}F''\left(\delta\right)+\ldots
\end{equation}
is also an asymptotic development of $R$. This question has also
been treated by Riemann; and again without the necessary residual
estimates. But since the series (46), because of its larger residual
terms, has no such theoretical and practical importance as the original
asymptotic developments, rather laborious investigation of the error
is omitted in the following presentation; perhaps this {[}presentation?{]}
brings out still more Riemann's formal power.

The formula (30), which can also be written in the form of
\[
\intop_{\mathclap{m\nwarrow m+1}}\frac{e^{-\pi i\left(x-\eta\right)^{2}+2\pi i\left(x-\eta\right)\left(\eta-m\right)}}{e^{2\pi ix}-1}dx=F\left(\delta\right)e^{-\frac{\pi i}{8}-\pi i\left(\eta-m\right)^{2}}
\]
allows the inversion
\begin{equation}
\frac{2}{\sqrt{2\pi}}e^{-\frac{\pi i}{8}-\pi i\left(\eta-m\right)^{2}}\intop_{\mathclap{0\swarrow1}}F\left(u+\delta\right)e^{iu^{2}-2\sqrt{2\pi}i\left(x-\eta\right)u}du
\end{equation}
\[
=\frac{e^{-\pi i\left(x-\eta\right)^{2}+2\pi i\left(x-\eta\right)\left(\eta-m\right)}}{e^{2\pi ix}-1},
\]
if $m<\mathfrak{R}\left(x\right)<m+1$. This results either by application
of Fourier's theorem or by transition to the complex conjugate variables
in (5). From (47) follows
\begin{equation}
\frac{e^{-2\pi imx}}{e^{2\pi ix}-1}=\left(-1\right)^{m}\frac{2}{\sqrt{2\pi}}e^{-\pi ix^{2}-\frac{\pi i}{8}}\intop_{\mathclap{0\swarrow1}}F\left(u+\delta\right)e^{i\left(u+\tau-\sqrt{2\pi}x\right)^{2}}du,
\end{equation}
also valid for $m<\mathfrak{R}\left(x\right)<m+1$. It is straightforward
to substitute for $F\left(u+\delta\right)$ herein the series 
\[
F\left(\delta\right)+\frac{F'\left(\delta\right)}{1!}u+\frac{F''\left(\delta\right)}{2!}u^{2}+\ldots
\]
and to calculate the contribution supplied by a single member of this
series by using the integral occurring in (16) 
\[
\intop_{\mathclap{\eta+\varepsilon\frac{\eta}{2}}}^{\mathclap{\eta-\varepsilon\frac{\eta}{2}}}g\left(x\right)dx=\intop_{\mathclap{\eta+\varepsilon\frac{\eta}{2}}}^{\mathclap{\eta-\varepsilon\frac{\eta}{2}}}x^{s-1}\frac{e^{-2\pi imx}}{e^{-2\pi ix}-1}dx
\]
\newpage\noindent In this way we find the asymptotic development
\begin{equation}
\intop_{C_{2}}g\left(x\right)dx
\end{equation}
\[
\sim\left(-1^{m}\right)\frac{2}{\sqrt{2\pi}}e^{-\frac{\pi i}{8}}\sum_{n=0}^{\infty}\frac{F^{\left(n\right)}\left(\delta\right)}{n!}\intop_{\mathclap{m\nwarrow m+1}}x^{s-1}e^{-\pi ix^{2}}\left\{ \;\intop_{\mathclap{0\swarrow1}}u^{n}e^{i\left(u+\tau-\sqrt{2\pi}x\right)^{2}}du\right\} dx.
\]
On the other hand, according to (13) and (44)
\begin{equation}
\intop_{C_{2}}g\left(x\right)dx=\left(-1\right)^{m}\left(\frac{t}{2\pi}\right)^{-\frac{1}{4}}e^{\vartheta i}R\left(1-e^{\pi is}\right).
\end{equation}
Now, since one can easily see that representation of $A_{n}$ as a
homogeneous linear function of $F^{\left(k\right)}\left(\delta\right)$
with constant coefficients as given in (41) is possible in only one
way, there follows from (49) and (50) the equation 
\begin{equation}
n!D_{n}\left(1-e^{\pi is}\right)
\end{equation}
\begin{flalign*}
 &  &  &  &  & =\frac{2}{\sqrt{2\pi}}\left(\frac{t}{2\pi}\right)^{\frac{1}{4}}e^{-\vartheta i-\frac{\pi i}{8}}\intop_{\mathclap{0\nwarrow1}}x^{s-1}e^{-\pi ix^{2}}\left\{ \;\intop_{\mathclap{0\swarrow1}}u^{n}e^{i\left(u+\tau-\sqrt{2\pi}x\right)^{2}}du\right\} dx & \left(n=0,1,2,\ldots\right),
\end{flalign*}
and in particular, if we set
\[
\frac{1}{\sqrt{2\pi}}\left(\frac{t}{2}\right)^{\frac{1}{4}}e^{\frac{\pi}{4}t}\sqrt{\varGamma\left(\frac{1}{4}+\frac{ti}{2}\right)\varGamma\left(\frac{1}{4}-\frac{ti}{2}\right)}=e^{\omega}
\]
then
\begin{equation}
\left\{ \begin{aligned}D_{0} & =e^{\omega}\vphantom{\frac{1}{1}}\\
D_{1} & =-\tau\left(e^{\omega}-e^{-\omega}\right)+\dfrac{\tau e^{\pi is-\omega}}{1-e^{\pi is}}\sim-\tau\left(e^{\omega}-e^{-\omega}\right).
\end{aligned}
\right.
\end{equation}
For the remaining $D_{n}$ a recurrence can be derived from (51) by
partial integration; but one may also obtain this in the following
way without additional calculation. Following (36), (37), (38)
\[
S\sim d_{0}F\left(\delta\right)+\frac{d_{1}}{1!}F'\left(\delta\right)+\frac{d_{2}}{2!}F''\left(\delta\right)+\cdots,
\]
where, by (38) and (39), the $d_{n}$ satisfy the recurrence
\begin{flalign*}
 &  &  &  &  & \tau d_{n}+\frac{1}{2}d_{n+1}+\frac{\left(n-1\right)\left(n-2\right)}{8}d_{n-3}=0 & \left(n=1,2,3,\ldots\right)
\end{flalign*}
Because
\[
e^{\left(\frac{t}{2}\log\frac{t}{2\pi}-\frac{t}{2}-\frac{\pi}{8}-\vartheta\right)i}S=R
\]
we have therefore for the $D_{n}$ the recursion formula
\begin{flalign}
 &  &  &  &  & D_{n+1}=-\frac{2}{n+1}\tau D_{n}-\frac{1}{4n\left(n+1\right)}D_{n-3} & \left(n=1,2,3,\ldots\right)
\end{flalign}
\newpage\noindent with $D_{-2}=0,D_{-1}=0$. From this one derives
with the help of (52) the values
\begin{align*}
D_{2} & =-\tau D_{1}\sim\tau^{2}\left(e^{\omega}-e^{-\omega}\right)\\
\\
D_{3} & =-\frac{2}{3}\tau D_{2}\sim-\frac{2}{3}\tau^{3}\left(e^{\omega}-e^{-\omega}\right)\\
\\
D_{4} & =-\frac{1}{2}\tau D_{3}-\frac{1}{2^{4}\cdot3}D_{0}\sim\frac{1}{3}\tau^{4}\left(e^{\omega}-e^{-\omega}\right)-\frac{1}{2^{4}\cdot3}e^{\omega}.
\end{align*}
One obtains the asymptotic developments of $D_{0},D_{1},\ldots$ themselves
from (43); according to this formula, specifically
\[
\omega\sim\frac{1}{8}\sum_{n=1}^{\infty}\frac{E_{n}}{n}\left(2t\right)^{-2n}=\frac{1}{2^{5}}t^{-2}+\frac{5}{2^{8}}t^{-4}+\frac{61}{2^{9}\cdot3}t^{-6}+\cdots.
\]
Substituting this in the above results for $D_{0},\ldots,D_{4}$,
it follows that
\begin{equation}
\left\{ \begin{aligned}D_{0} & \sim1+\tfrac{1}{2^{5}}\tau^{-4}+\tfrac{41}{2^{11}}\tau^{-8}+\cdots\\
\\
D_{1} & \sim-\tfrac{1}{2^{4}}\tau^{-3}-\tfrac{5}{2^{7}}\tau^{-7}+\cdots\\
\\
D_{2} & \sim\tfrac{1}{2^{4}}\tau^{-2}+\tfrac{5}{2^{7}}\tau^{-6}+\cdots\\
\\
D_{3} & \sim-\tfrac{1}{2^{3}\cdot3}\tau^{-1}-\tfrac{5}{2^{6}\cdot3}\tau^{-5}+\cdots\\
\\
D_{4} & \sim\tfrac{19}{2^{9}\cdot3}\tau^{-4}+\cdots.
\end{aligned}
\right.
\end{equation}
From the recurrence relation (53) it follows that all power exponents
occurring in the asymptotic series for $D_{n}$ are $\equiv n\,\left(\mathrm{mod}\,4\right)$.
Accordingly, the order of all derivatives of $F\left(\delta\right)$
occurring in $C_{n}$ are of the form $3n\lyxmathsym{\textendash}4k$,
as can easily be confirmed in the expressions found for $C_{0},C_{1},C_{2},C_{3},C_{4}$.
If one writes
\[
R\sim\sum b_{kl}\,F^{\left(3l-4k\right)}\left(\delta\right)\tau^{-l},
\]
where the summation index $k$ passes through the values $0,\ldots,\left[\frac{3l}{4}\right]$
and the summation index $l$ the values $0,1,\ldots$, then all $b_{kl}$
with $l\leq4$ are determined, while, according to (54) the values
$b_{00},b_{34},b_{68},b_{23},b_{57},b_{12},b_{46},b_{01},b_{35},b_{24}$
are known; one sees immediately that the values for $b_{00},b_{34},b_{23},b_{12},b_{01},b_{24}$
occurring in both (45) and (54) agree.

For the numerical computation of the $b_{kl}$, and for practical
applications of asymptotic series, ordering by increasing powers of
$\tau^{-1}$ is preferred. The determination of the $D_{n}$ by (53)
is in fact more laborious than the previously treated determination
of the $C_{n}$; moreover, the consecutive $D_{n}$ do not have monotonically
decreasing magnitudes, but $D_{3n-2},D_{3n-1},D_{3n}$ do have the
exact

\newpage\noindent orders $\tau^{-\left(n+2\right)},\tau^{-\left(n+1\right)},\tau^{-n}$,
so that one must therefore, for example, obtain $D_{5}$ to $D_{12}$
only up to the previous error $O\left(\tau^{-5}\right)$.

The transition to the $D_{n}$ was done using the formula (48). If
one tried to obtain from (48) an exact expression for $\zeta\left(s\right)$
and not merely an asymptotic series, then one would come to the approach
that will be discussed in the next section.

\section{The Integral Representation of the Zeta Function.}

The explicit determination of the coefficients in the asymptotic series
for $\zeta\left(s\right)$ was based on the formula (5) of § 1. With
the help of this formula Riemann has derived a further rather interesting
expression for $\zeta\left(s\right)$, which appears to have escaped
the attention of other mathematicians until 1926.

For now, let $\sigma<0$, and let $u^{-s}$ take the principal value
in the $u$ plane cut from 0 to $-\infty$. One multiplies (5) by
$u^{-s}$ and integrates over $u$ from 0 to $e^{\frac{\pi i}{4}\infty}$
along the bisector of the first quadrant. Now, if the abbreviation
$e^{\frac{\pi i}{4}}=\bar{\varepsilon}$ is set,
\begin{align*}
\intop_{0}^{\bar{\varepsilon}\infty}\frac{u^{-s}}{1-e^{-2\pi iu}}du & =-\intop_{0}^{\bar{\varepsilon}\infty}u^{-s}\sum_{n=1}^{\infty}e^{2\pi inu}du=-\sum_{n=1}^{\infty}\intop_{0}^{\bar{\varepsilon}\infty}u^{-s}e^{2\pi inu}du\\
\\
 & =-\varGamma\left(1-s\right)\sum_{n=1}^{\infty}\left(2\pi ne^{-\frac{\pi i}{2}}\right)^{s-1}\\
\\
 & =-\left(2\pi\right)^{s-1}e^{\frac{\pi i}{2}\left(1-s\right)}\varGamma\left(1-s\right)\zeta\left(1-s\right)
\end{align*}
and
\begin{align*}
\intop_{0}^{\bar{\varepsilon}\infty}u^{-s}\left(\,\;\intop_{\mathclap{0\nwarrow1}}\frac{e^{\pi ix^{2}+2\pi iux}}{e^{\pi ix}-e^{-\pi ix}}dx\right)du & =\intop_{\mathclap{0\nwarrow1}}\frac{e^{-\pi ix^{2}}}{e^{\pi ix}-e^{-\pi ix}}\left(\intop_{0}^{\bar{\varepsilon}\infty}u^{-s}e^{2\pi iux}du\right)dx\\
 & =\left(2\pi\right)^{s-1}e^{\frac{\pi i}{2}\left(1-s\right)}\varGamma\left(1-s\right)\intop_{\mathclap{0\nwarrow1}}\frac{e^{-\pi ix^{2}}x^{s-1}}{e^{\pi ix}-e^{-\pi ix}}dx,
\end{align*}
so by (5)
\begin{align}
\left(2\pi\right)^{s-1}e^{\frac{\pi i}{2}\left(1-s\right)}\varGamma\left(1-s\right)\left\{ \zeta\left(1-s\right)\vphantom{+\intop_{\mathclap{0\nwarrow1}}\frac{x^{s-1}e^{-\pi ix^{2}}}{e^{\pi ix}-e^{-\pi ix}}dx}\right. & \left.+\intop_{\mathclap{0\nwarrow1}}\frac{x^{s-1}e^{-\pi ix^{2}}}{e^{\pi ix}-e^{-\pi ix}}dx\right\} \\
 & +\intop_{0}^{\bar{\varepsilon}\infty}\frac{u^{-s}e^{\pi iu^{2}}}{e^{\pi iu}-e^{-\pi iu}}du=0.\nonumber 
\end{align}
\newpage\noindent Here the second integral may be put in the form
\[
\frac{1}{e^{\pi is}-1}\intop_{\mathclap{0\swarrow1}}\frac{u^{-s}e^{\pi iu^{2}}}{e^{\pi iu}-e^{-\pi iu}}du
\]
where by the symbol $0\swarrow1$ denotes the integration path that
arises from the first integral by reflection in the real axis. Multiplying
(55) by the factor 
\[
2^{1-s}\pi^{\frac{1-s}{2}}e^{\frac{\pi i}{2}\left(s-1\right)}\frac{\varGamma\left(\frac{1-s}{2}\right)}{\varGamma\left(1-s\right)}
\]
and considering the relation
\[
\frac{2^{-s}\pi^{\frac{1-s}{2}}\varGamma\left(\frac{1-s}{2}\right)}{\sin\frac{\pi s}{2}\varGamma\left(1-s\right)}=\pi^{-\frac{s}{2}}\varGamma\left(\frac{s}{2}\right),
\]
one obtains, valid throughout the $s$-plane, the formula 
\begin{align}
\pi^{-\frac{1-s}{2}}\varGamma\left(\frac{1-s}{2}\right)\zeta\left(1-s\right) & =\pi^{-\frac{s}{2}}\varGamma\left(\frac{s}{2}\right)\intop_{\mathclap{0\swarrow1}}\frac{x^{-s}e^{\pi ix^{2}}}{e^{\pi ix}-e^{-\pi ix}}dx\\
 & +\pi^{-\frac{1-s}{2}}\varGamma\left(\frac{1-s}{2}\right)\intop_{\mathclap{0\searrow1}}\frac{x^{s-1}e^{-\pi ix^{2}}}{e^{\pi ix}-e^{-\pi ix}}dx.\nonumber 
\end{align}
Riemann does not write everything in this symmetrical form; but the
version chosen here seems to be expedient for applications. One has
now made the functional equation for $\zeta(s)$ evident; because
for $\sigma=\frac{1}{2}$  the two terms on the right side are complex
conjugates, so $\pi^{-\frac{s}{2}}\varGamma\left(\frac{s}{2}\right)\zeta\left(s\right)$
is real there, and since this function is real for $\sigma>1$, so
by the principle of symmetry, the functional equation 
\begin{equation}
\pi^{-\frac{s}{2}}\varGamma\left(\frac{s}{2}\right)\zeta\left(s\right)=\pi^{-\frac{1-s}{2}}\varGamma\left(\frac{1-s}{2}\right)\zeta\left(1-s\right)
\end{equation}
applies for $\sigma=\frac{1}{2}$  and thus generally for any $s$.

If one sets further
\begin{equation}
\mathfrak{f}\left(s\right)=\intop_{\mathclap{0\swarrow1}}\frac{x^{-s}e^{\pi ix^{2}}}{e^{\pi ix}-e^{-\pi ix}}dx
\end{equation}
\begin{equation}
\varphi\left(s\right)=2\pi^{-\frac{s}{2}}\varGamma\left(\frac{s}{2}\right)\mathfrak{f}\left(s\right)
\end{equation}
then by (56) and (57)
\begin{flalign}
 &  &  &  &  & \pi^{-\frac{s}{2}}\varGamma\left(\frac{s}{2}\right)\zeta\left(s\right)=\mathfrak{R}\left(\varphi\left(s\right)\right)\qquad & \left(\sigma=\tfrac{1}{2}\right);
\end{flalign}
\newpage\noindent with this the study of $\zeta\left(s\right)$
on the critical line is reduced to the study of the real part of $\varphi\left(s\right)$.

\section{The Significance of the Two Riemann Formulas for the Theory of the
Zeta Function.}

The principal term of the asymptotic series for $\zeta\left(s\right)$,
based on equation (32), is the expression
\begin{flalign*}
 &  &  &  &  &  &  & \sum_{l=1}^{m}l^{-s}+\frac{\left(2\pi\right)^{s}}{\pi}\sin\frac{\pi s}{2}\varGamma\left(1-s\right)\sum_{l=1}^{m}l^{s-1} & \left(m=\left[\sqrt{\frac{t}{2\pi}}\;\,\right]\right),
\end{flalign*}
which was also found by Hardy and Littlewood, but in place of Riemann's
development for $S$, they have only given an upper bound for its
absolute value. They also discovered a more general form of the principal
term, namely
\begin{equation}
\sum_{l\leq x}l^{-s}+\frac{\left(2\pi\right)^{s}}{\pi}\sin\frac{\pi s}{2}\varGamma\left(1-s\right)\sum_{l\leq y}l^{s-1}
\end{equation}
with $xy=\frac{t}{2\pi}$. This was not given previously by Riemann;
but one can imagine that, without any major difficulty, one can arrive
along Riemann's lines at a full asymptotic development of the expression
(61), which plays the same role for the function $\varPhi\left(\tau,u\right)$,
defined by (7), as the special function $\varPhi\left(-1,u\right)$
of Riemann.

For the applications for which Hardy and Littlewood have made use
of their formula, especially for estimating the number of zeros $N_{0}\left(T\right)$
of $\zeta\left(\frac{1}{2}+ti\right)$  lying in the interval $0<t<T$,
for which the Riemann formula gives more accurate values, it seems
not to be a better result. In the above-mentioned place, however,
Riemann has claimed that $N_{0}\left(T\right)$ is asymptotically
equal to $\frac{T}{2\pi}\log\frac{T}{2\pi}-\frac{T}{2\pi}$, and so
asymptotically equal to the number of all zeros $N\left(T\right)$
of $\zeta\left(s\right)$ lying in the strip $0<t<T$, and that this
can be proven with the help of his new development; but from his Nachlass
it is not clear whether he had devised this proof. In the representation,
valid on $\sigma=\frac{1}{2}$ 
\begin{equation}
e^{\vartheta i}\zeta\left(\frac{1}{2}+ti\right)=2\sum_{n=1}^{m}\frac{\cos\left(\vartheta-t\log n\right)}{\sqrt{n}}+O\left(t^{-\frac{1}{4}}\right)
\end{equation}
\[
\vartheta=\frac{t}{2}\log\frac{t}{2\pi}-\frac{t}{2}+O\left(1\right)
\]
\newpage\noindent the first term of the right-side trigonometric
sum, namely $\cos\vartheta$, actually has asymptotically $\frac{T}{2\pi}\log\frac{T}{2\pi}-\frac{T}{2\pi}$
zeros in the interval $0<t<T$; and the coefficients $\frac{1}{\sqrt{1}},\frac{1}{\sqrt{2}},\frac{1}{\sqrt{3}},\ldots$
decrease monotonically. Perhaps Riemann believed that this observation
could be used in the proof of his assertion.

It is obvious to use the exact Riemann formula for estimating the
mean values
\begin{flalign*}
 &  &  &  &  & \frac{1}{T}\intop_{0}^{T}\left|\zeta\left(\frac{1}{2}+ti\right)\right|^{2n}dt & \left(n=3,4,\ldots\right);
\end{flalign*}
these averages are well known and are closely related to the so-called
Lindelöf conjecture. But here one encounters considerable difficulties
of an arithmetic nature, resulting from the factorization of natural
numbers.

For the establishment of a numerical table of the zeta function, in
particular for the further calculation of the zeros, the asymptotic
development is of great benefit. However, if one would have to use
it for the purposes of practical application, a more careful estimate
of the remainder term must be developed, than is derived in § 2. Riemann
employed his formula for rather extensive calculations for determining
positive zeros of $\zeta\left(\frac{1}{2}+ti\right)$. For the smallest
positive zero he finds the value of $a_{1}=$ 14.1386; Gram has calculated
the value 14.1347, different by less than 3 parts per thousand. A
lower bound for \foreignlanguage{ngerman}{$a_{1}$} also results using
the product representation of $\zeta\left(s\right)$, yielding the
easily proven equation
\[
\sum_{n=1}^{\infty}\left(a_{n}^{2}+\dfrac{1}{4}\right)^{-1}=1+\dfrac{1}{2}\gamma-\dfrac{1}{2}\log\pi-\log2,
\]
where $\gamma$ is Euler's constant, and where $a_{n}$ passes through
all solutions $a$ of $\zeta\left(\frac{1}{2}+ai\right)=0$ located
in the right half-plane. From this Riemann asserts
\[
\sum_{n=1}^{\infty}\left(a_{n}^{2}+\dfrac{1}{4}\right)^{-1}=0.02309\:57089\:66121\:03381.
\]
For $a_{3}$ he finds the value 25.31; while Gram gives 25.01 for
this.

The second Riemann formula, namely the integral representation of
$\zeta\left(s\right)$, may perhaps be of greater interest for the
theory. One will try to get, from (60), some information about the
distribution\newpage\noindent of the zeros of $\zeta\left(s\right)$
on the critical line. Let $t$ increase through an interval $t_{1}<t<t_{2}$.
During this process $\arg\,\varphi\left(\frac{1}{2}+ti\right)$ increases
by $A$, where the change when passing a possible zero of $\varphi\left(s\right)$
located on $\sigma=\frac{1}{2}$ equals $\pi$ multiplied by the multiplicity
of this zero, so by (60) the number of zeros of $\zeta\left(\frac{1}{2}+ti\right)$
within $t_{1}<t<t_{2}$ is greater than $\frac{\left|A\right|}{\pi}-1$.
But now, however, we have 
\begin{equation}
\arg\left\{ \pi^{-\frac{s}{2}}\varGamma\left(\frac{s}{2}\right)\right\} =\vartheta=\dfrac{t}{2}\log\dfrac{t}{2\pi}-\dfrac{t}{2}+O\left(1\right),
\end{equation}
and thus according to (59) the number of zeros of $\zeta\left(\frac{1}{2}+ti\right)$
lying in the interval $0<t<T$ is asymptotically at least equal to
$\frac{T}{2\pi}\log T$, ie. asymptotically equal to the number zeros
of $\zeta\left(s\right)$ lying in the strip $0<t<T$, if the argument
of the function $\mathfrak{f}\left(\frac{1}{2}+ti\right)$ defined
by (58) decreases more slowly for $t\rightarrow\infty$ than $-t\log t$.
For each half-strip $\sigma_{1}\leq\sigma\leq\sigma_{2},t>0$ one
can develop $\mathfrak{f}\left(s\right)$ into an asymptotic series by
the method of § 2; but one again obtains as a principal term a sum
of $\left[\sqrt{\frac{t}{2\pi}}\;\,\right]$ summands, that is $\sum_{n=1}^{m}n^{-s}$;
and the study of the argument of this sum is a problem of exactly
the same difficulty as the study of the zeros of the sum occurring
in (62), so that therefore with the introduction of $\mathfrak{f}\left(s\right)$ we
have apparently gained nothing.

If one now considers the rectangle with the sides $\sigma=\frac{1}{2},\sigma=2,t=0,t=T$, of
which the upper side contains no zeros of $\mathfrak{f}\left(s\right)$,
then the change of $\frac{1}{2\pi}\arg\,\mathfrak{f}\left(s\right)$
for positive circulation of the rectangle equals the number of zeros
of $\mathfrak{f}\left(s\right)$ within the rectangle. On the lower
side $\arg\,\mathfrak{f}\left(s\right)$ changes by $O\left(1\right)$ and
on the right side, which obeys the asymptotic series, also only by
$O\left(1\right)$. Furthermore, it can be shown by the usual methods
in the theory of the zeta function that the change on the upper side
is at most $O\left(\log T\right)$. Consequently, except for an error
of the order of $\log T$, the change of $\arg\,\mathfrak{f}\left(\frac{1}{2}+ti\right)$  in
the interval $0<t<T$ is equal to $-2\pi$ multiplied by the number
of zeros of $\mathfrak{f}\left(s\right)$ within the rectangle. With
this, the problem is reduced to the study of the zeros of the entire
transcendent $\mathfrak{f}\left(s\right)$.\newpage Riemann tried
to obtain a statement concerning the zeros of $\mathfrak{f}\left(s\right)$,
so to that end he formed from (58) the relation
\[
\left|\mathfrak{f}\left(\sigma+ti\right)\right|^{2}=\intop_{\mathclap{0\swarrow1}}\;\;\;\intop_{\mathclap{0\searrow1}}\dfrac{x^{-\sigma-ti}y^{-\sigma+ti}e^{\pi i\left(x^{2}-y^{2}\right)}}{\left(e^{\pi ix}-e^{-\pi ix}\right)\left(e^{\pi iy}-e^{-\pi iy}\right)}dxdy
\]
and brought the complex double integral into a different form by introduction
of new variables, deformation of the integration area and use of the
residue theorem; however, this lead to no useful result.

So far very little is known about the locations of the zeros of $\mathfrak{f}\left(s\right)$.
By Riemann we find no further remarks on this subject; in the context
of this historical-mathematical treatise, therefore, the following
remarks on the theory of $\mathfrak{f}\left(s\right)$ must be kept
brief. They provide a proof of the inequality
\[
N_{0}\left(T\right)>\dfrac{3}{8\pi}e^{-\frac{3}{2}}T+o\left(T\right).
\]
For $\mathfrak{f}\left(s\right)$ one can find an asymptotic series
by the method of § 2; for the present purpose it is sufficient to
consider the principal term of this series. First it will be shown
that, in regions where $t>0,-\sigma\geq t^{\frac{3}{7}}$, the following
formula applies
\begin{flalign}
 &  &  &  &  & \mathfrak{f}\left(s\right)\sim e^{\frac{\pi i}{4}\left(s-\frac{7}{2}\right)}\pi^{\frac{s-1}{2}}\sin\dfrac{\pi s}{2}\varGamma\left(\frac{1-s}{2}\right)\dfrac{\sin\pi\eta}{\cos2\pi\eta} & \left(\left|s\right|\rightarrow\infty\right)
\end{flalign}
where the abbreviation
\begin{flalign*}
 &  &  &  &  & \eta=\sqrt{\dfrac{s-1}{2\pi i}} & \left(0<\arg\,\eta<\tfrac{\pi}{4}\right)
\end{flalign*}
is used.

Now, from (56)
\begin{equation}
\mathfrak{f}\left(s\right)=\pi^{s-\frac{1}{2}}\dfrac{\varGamma\left(\dfrac{1-s}{2}\right)}{\varGamma\left(\dfrac{s}{2}\right)}\left\{ \zeta\left(1-s\right)-\intop_{\mathclap{0\searrow1}}\dfrac{x^{s-1}e^{-\pi ix^{2}}}{e^{\pi ix}-e^{-\pi ix}}dx\right\} .
\end{equation}
The saddle point of the function $x^{s-1}e^{-\pi ix^{2}}$ is located
at $x=\eta$.   One sets
\[
\mathfrak{R}\left(\eta\right)=\eta_{1},\qquad\mathfrak{J\left(\eta\right)=\eta_{2}}
\]
\[
m=\left[\eta_{1}+\eta_{2}\right]
\]
\[
z=x-\eta
\]
\[
w\left(z\right)=e^{2\pi i\eta^{2}\left\{ \log\left(1+\frac{z}{\eta}\right)-\frac{z}{\eta}+\frac{1}{2}\left(\frac{z}{\eta}\right)^{2}\right\} }-1.
\]
\newpage\noindent For every natural number $k$ it follows from
Cauchy's theorem 
\begin{equation}
\intop_{\mathclap{0\searrow1}}\dfrac{x^{s-1}e^{-\pi ix^{2}}}{e^{\pi ix}-e^{-\pi ix}}dx
\end{equation}
\[
={\displaystyle \sum_{n=1}^{k}}n^{s-1}+\eta^{s-1}e^{-\pi i\eta^{2}}\left\{ \quad\intop_{\mathclap{k\searrow k+1}}\dfrac{e^{-2\pi i\left(x-\eta\right)^{2}}}{e^{\pi ix}-e^{-\pi ix}}dx+\intop_{\mathclap{k\searrow k+1}}\dfrac{e^{-2\pi i\left(x-\eta\right)^{2}}}{e^{\pi ix}-e^{-\pi ix}}w\left(z\right)dx\right\} .
\]
One is now tempted to proceed exactly by the method of § 2, and so
one would choose $k=m$; but then one would obtain (64) directly only
in the smaller rectangular region $t>0,\;\lyxmathsym{\textendash}\sigma\geq t^{\frac{1}{2}}$,
and the extension to the additional region $t^{\frac{1}{2}}>\lyxmathsym{\textendash}\sigma\geq t^{\frac{3}{7}}$ requires
the elimination of certain additional terms. That is why we initially
let $k$ be chosen arbitrarily.

The first integral on the right side of (66) can be calculated according
to Riemann's method of § 1; one obtains
\begin{equation}
\intop_{\mathclap{k\searrow k+1}}\dfrac{e^{-2\pi i\left(x-\eta\right)^{2}}}{e^{\pi ix}-e^{-\pi ix}}dx=\dfrac{\sqrt{2}e^{\frac{3\pi i}{8}}\sin\pi\eta+\left(-1\right)^{k-1}e^{2\pi i\eta-2\pi i\left(\eta-k\right)^{2}}}{2\cos2\pi\eta}.
\end{equation}
In the second integral we lead the integration path through the saddle
point $x=\eta$ and traverse it parallel to the bisector of the second
and fourth quadrants. Thus it crosses the real axis at the point $\eta_{1}+\eta_{2}$.
However, in order to avoid the neighborhoods of the Poles at $x=m~~$
and $x=m+1~~$, one can still replace the parts of the path of integration
lying inside of the circles $\left|x-m\right|=\frac{1}{2}$ and $\left|x-m-1\right|=\frac{1}{2}$ by
arcs of these circles. If one makes the assumption that
\[
k=m+r\geq m,
\]
then 
\begin{align}
\intop_{\mathclap{k\searrow k+1}}\dfrac{e^{-2\pi i\left(x-\eta\right)^{2}}}{e^{\pi ix}-e^{-\pi ix}}w\left(z\right)dx & =\sum_{l=1}^{r}\left(-1\right)^{m+l-1}e^{-2\pi i\left(m+l-\eta\right)^{2}}w\left(m+l-\eta\right)\\
 & +\intop_{\mathclap{m\searrow m+1}}\dfrac{e^{-2\pi i\left(x-\eta\right)^{2}}}{e^{\pi ix}-e^{-\pi ix}}w\left(z\right)dx.\nonumber 
\end{align}
For $w\left(z\right)$ one needs two estimates. The first refers to
the circle $\left|z\right|\leq\frac{1}{2}\left|\eta\right|$; specifically,
this is
\[
\left|\log\left(1+\frac{z}{\eta}\right)-\frac{z}{\eta}+\frac{1}{2}\left(\frac{z}{\eta}\right)^{2}\right|
\]
\[
=\left|\sum_{n=3}^{\infty}\dfrac{\left(-1\right)^{n-1}}{n}\left(\frac{z}{\eta}\right)^{n}\right|\leq\dfrac{1}{3}\left|\frac{z}{\eta}\right|^{3}\dfrac{1}{1-\left|\frac{z}{\eta}\right|}\leq\dfrac{2}{3}\left|\frac{z}{\eta}\right|^{3}
\]
 \newpage\noindent and hence 
\begin{flalign}
 &  &  &  &  & \left|w\left(z\right)\right|\leq e^{\frac{4\pi}{3}\left|\frac{z^{3}}{\eta}\right|}-1\qquad & \left(\left|z\right|\leq\tfrac{1}{2}\left|\eta\right|\right).
\end{flalign}
The second refers to the parts of the integration path located outside
this circle. If one now sets $\mathfrak{R\left(\mathit{ze^{\frac{\pi i}{4}}}\right)=\mathit{u}}$,
$\mathfrak{I}\left(\mathit{ze^{\frac{\pi i}{4}}}\right)=\mathit{v}$,
then on the path of integration $-\frac{1}{2}\leq v\leq+\frac{1}{2}$,
and in the case of $\left|\eta\right|>1$ there applies outside the
circle $\left|z\right|=\frac{1}{2}\left|\eta\right|$ the inequality
\[
\left|\frac{v}{u}\right|<\left(\left|\eta\right|^{2}-1\right)^{-\frac{1}{2}},
\]
so
\begin{flalign*}
 &  &  & \qquad\arg\left(1+\dfrac{iv}{u}\right)\rightarrow0 & \left(\left|s\right|\rightarrow\infty\right)
\end{flalign*}
and
\[
\dfrac{\pi}{4}-\varepsilon<\left|\arg\,\dfrac{z}{\eta}\right|<\dfrac{3\pi}{4}+\varepsilon,
\]
with $\varepsilon\rightarrow0$ for $\left|s\right|\rightarrow\infty$.
But then
\[
\left|2\pi i\eta^{2}\left\{ \log\left(1+\dfrac{z}{\eta}\right)-\dfrac{z}{\eta}+\dfrac{1}{2}\left(\dfrac{z}{\eta}\right)^{2}\right\} \right|=2\pi\left|\eta\right|^{2}\cdot\left|{\displaystyle \intop_{0}^{\frac{z}{\eta}}\dfrac{x^{2}}{1+x}dx}\right|
\]
\[
\leq2\pi\left|\eta\right|^{2}\intop_{0}^{\left|\frac{z}{\eta}\right|}\dfrac{x}{\sin\left(\dfrac{\pi}{4}-\varepsilon\right)}dx=\dfrac{\pi\left|z\right|^{2}}{\sin\left(\dfrac{\pi}{4}-\varepsilon\right)}\leq\dfrac{3}{2}\pi\left|z\right|^{2}
\]
and\foreignlanguage{ngerman}{
\[
\left|w\left(z\right)\right|<e^{\frac{3}{2}\pi\left|z\right|^{2}}.
\]
}Further, on the integration path
\[
\left|e^{-2\pi iz^{2}}\right|=e^{-2\pi\left(u^{2}-v^{2}\right)}\leq e^{\pi-2\pi\left|z\right|^{2}};
\]
and therefore
\[
\qquad\qquad\qquad\qquad\qquad\qquad\left|e^{-2\pi iz^{2}}w\left(z\right)\right|\leq\left\{ \begin{aligned} & e^{\pi-\frac{\pi}{2}\left|z\right|^{2}} &  &  &  &  &  &  &  &  &  &  &  &  & \left(\left|z\right|>\tfrac{1}{2}\left|\eta\right|\right)\\
 & e^{\pi-2\pi\left|z\right|^{2}}\left(e^{\frac{4\pi}{3}\left|\frac{z^{3}}{\eta}\right|}-1\right) &  &  &  &  &  &  &  &  &  &  &  &  & \left(\left|z\right|\leq\tfrac{1}{2}\left|\eta\right|\right).
\end{aligned}
\right.
\]
This yields
\[
{\displaystyle \intop_{\mathclap{m\searrow m+1}}\dfrac{e^{-2\pi i\left(x-\eta\right)^{2}}}{e^{\pi ix}-e^{-\pi ix}}w\left(z\right)dx=O\left(e^{-\pi\eta_{_{2}}}\,\eta^{-1}\right)}
\]
and together with (65), (66), (67), (68)\newpage{} 
\begin{equation}
\mathfrak{f}\left(s\right)=\pi^{s-\frac{1}{2}}\dfrac{\varGamma\left(\dfrac{1-s}{2}\right)}{\varGamma\left(\dfrac{s}{2}\right)}\eta^{s-1}e^{-\pi i\eta^{2}}\left\{ e^{\pi i\eta^{2}}{\displaystyle \sum_{n=m+r+1}^{\infty}\left(\dfrac{n}{\eta}\right)^{s-1}}\right.
\end{equation}
\[
-\dfrac{\sqrt{2}e^{\frac{3\pi i}{8}}\sin\pi\eta+\left(-1\right)^{m+r-1}e^{2\pi i\eta-2\pi i\left(\eta-m-r\right)^{2}}}{2\cos2\pi\eta}
\]
\[
\left.+\sum_{l=1}^{r}\left(-1\right)^{m+l}e^{-2\pi i\left(m+l-\eta\right)^{2}}w\left(m+l-\eta\right)+O\left(e^{-\pi\eta_{_{2}}}\,\eta^{-1}\right)\right\} \,.
\]

One has now to show, with a suitable choice of $r$ and for $\left|s\right|\rightarrow\infty$
in the region $t>0$, $-\sigma\geq t^{\frac{3}{7}}$, that the term
\[
-\dfrac{\sqrt{2}e^{\frac{3\pi i}{8}}\sin\pi\eta}{2\cos2\pi\eta}
\]
is of higher order than the other terms in the braces. Firstly,
\begin{align}
\left|e^{\pi i\eta^{2}}{\displaystyle \sum_{n=m+r+1}^{\infty}\left(\dfrac{n}{\eta}\right)^{s-1}}\right| & <e^{-2\pi\eta_{1}\eta_{2}}\left|\eta^{1-s}\right|\left\{ \dfrac{\left(m+r+1\right)^{\sigma}}{-\sigma}+\left(m+r+1\right)^{\sigma-1}\right\} \\
 & <e^{-2\pi\eta_{1}\eta_{2}+t\arg\eta}\left(\dfrac{m+r+1}{\left|\eta\right|}\right)^{\sigma-1}\left(\dfrac{m+r+1}{-\sigma}+1\right)\,;\nonumber 
\end{align}
for 
\begin{equation}
-2\pi\eta_{1}\eta_{2}+t\arg\eta<-2\pi\eta_{1}\eta_{2}+t\dfrac{\eta_{2}}{\eta_{1}}=-2\pi\dfrac{\eta_{2}^{3}}{\eta_{1}}<0
\end{equation}
and
\[
\left(\dfrac{m+1}{\left|\eta\right|}\right)^{\sigma-1}<\left(\dfrac{\eta_{1}^{2}+\eta_{2}^{2}+2\eta_{1}\eta_{2}}{\eta_{1}^{2}+\eta_{2}^{2}}\right)^{\frac{\sigma-1}{2}}<e{}^{\frac{\eta_{1}\eta_{2}}{\eta_{1}^{2}+\eta_{2}^{2}}\frac{\sigma-1}{2}}<e{}^{\frac{\eta_{2}}{2\eta_{1}}\cdot\frac{\sigma-1}{2}}=e^{-\pi\eta_{2}^{2}}
\]
therefore
\begin{equation}
\left|e^{\pi i\eta^{2}}{\displaystyle \sum_{n=m+1}^{\infty}\left(\dfrac{n}{\eta}\right)^{s-1}}\right|=O\left(e^{-\pi\eta_{2}^{2}}\left(1+\frac{\eta_{1}}{-\sigma}\right)\right);
\end{equation}
and also
\begin{equation}
\left|\left(-1\right)^{m-1}e^{2\pi i\eta-2\pi i\left(\eta-m\right)^{2}}\right|=e^{-2\pi\eta_{2}-4\pi\left(m-\eta_{1}\right)\eta_{2}}<e^{-4\pi\left(\eta_{2}-\frac{1}{2}\right)\eta_{2}}.
\end{equation}
Now, for the sub-region $t>0$, $-\sigma\geq t^{\frac{5}{8}}$ the
inequality
\[
\eta_{2}=\dfrac{1-\sigma}{4\pi\eta_{1}}>\dfrac{1-\sigma}{2}\left\{ 2\pi\left(t+1-\sigma\right)\right\} ^{-\frac{1}{2}}>\dfrac{1}{2}t^{\frac{5}{8}}\left\{ 2\pi\left(t+t^{\frac{5}{8}}\right)\right\} ^{-\frac{1}{2}}
\]
is fulfilled and the right side becomes infinite with $t$, so it
follows in view of (73) and (74) that the expression in the braces
of (70) in the aforementioned sub-region for $r=0$ takes the value\newpage
\begin{flalign}
 &  &  &  &  & \dfrac{\sqrt{2}e^{\frac{3\pi i}{8}}\sin\pi\eta}{2\cos2\pi\eta}\left(1+o\left(1\right)\right)\qquad & \left(\left|s\right|\rightarrow\infty\right)
\end{flalign}
 In the sub-region $t>0$, $t^{\frac{5}{8}}>-\sigma\geq t^{\frac{3}{7}}$,
now to be treated, one chooses
\[
r=\left[\left|\sigma\right|^{\frac{1}{5}}\right].
\]
Then for sufficiently large $t$
\[
\left(\dfrac{m+r+1}{\left|\eta\right|}\right)^{\sigma-1}<\left(\dfrac{\left|\eta\right|+r}{\left|\eta\right|}\right)^{\sigma-1}<e^{\frac{\sigma-1}{2}\cdot\frac{r}{2\eta_{1}}}=e^{-\pi r\eta_{2}},
\]
thus, by (71) and (72)
\begin{equation}
e^{\pi r\eta^{2}}{\displaystyle \sum_{n=m+r+1}^{\infty}}\left(\dfrac{n}{\eta}\right)^{s-1}=O\left(e^{-\pi r\eta_{2}}\left(1+\left|\dfrac{\eta}{\sigma}\right|\right)\right)=O\left(e^{-\frac{1}{2}t^{\frac{1}{70}}}\right).
\end{equation}
Moreover, for $l=1,\ldots,r$
\begin{equation}
\left|m+l-\eta\right|^{2}\leq\left(r+\eta_{2}\right)^{2}+\eta_{2}^{2}=O\left(\left|\sigma\right|^{\frac{2}{5}}\right)=O\left(t^{\frac{1}{4}}\right),
\end{equation}
so therefore $m+l-\eta$ for sufficiently large $t$ lies in the circle
$\left|z\right|\leq\frac{1}{2}\left|\eta\right|$ and (69) applies
for $z=m+l-\eta$; because of (77) it follows that
\begin{equation}
w\left(m+l-\eta\right)=O\left(\left|\sigma\right|^{\frac{3}{5}}\left|\eta\right|^{-1}\right).
\end{equation}
Finally, for $l=1,\:\ldots\:,\:r$
\begin{equation}
\left|e^{-2\pi i\left(m+l-\eta\right)^{2}}\right|<e^{-4\pi\left(\eta_{2}+l-1\right)\eta_{2}}\leq e^{-4\pi r\eta_{2}^{2}}
\end{equation}
and for sufficiently large $t$ 
\begin{equation}
\left|e^{-2\pi i\eta-2\pi i\left(m+l-\eta\right)^{2}}\right|<e^{-3\pi r\eta_{2}}=O\left(e^{-t^{\frac{1}{70}}}\right),
\end{equation}
thus, according to (78) and (79)
\begin{equation}
\sum_{l=1}^{r}\left(-1\right)^{m+l}e^{-2\pi i\left(m+l-\eta\right)^{2}}w\left(m+l-\eta\right)
\end{equation}
\[
=r\,O\left(e^{-4\pi\eta_{2}^{2}}\left|\sigma\right|^{\frac{3}{5}}\left|\eta\right|^{-1}\right)=O\left(e^{-4\pi\eta_{2}^{2}}\left|\sigma\right|^{\frac{4}{5}}\left|\eta\right|^{-1}\right).
\]
Considering now the inequalities
\begin{align*}
 & \left|\sin\pi\eta\right|\geq\sinh\pi\eta_{2}>\pi\eta_{2}>\dfrac{\left|\sigma\right|}{4\left|\eta\right|}\\
 & \left|\cos2\pi\eta\right|\leq\cosh2\pi\eta_{2}<2e^{2\pi\eta_{2}},
\end{align*}
\newpage\noindent the estimates (76), (80), (81) show that even
in the region $t>0$, $t^{\frac{5}{8}}>-\sigma\geq t^{\frac{3}{7}}$
the value in the braces of (70) is given by the expression (75).

The assertion in (64) then follows by applying Stirling's formula.

One can, incidentally, prove (64) even for the larger region $t>0$,
$-\sigma\geq t^{\varepsilon}$ where $\varepsilon$ is any fixed positive
number; but for the following it is sufficient for each value of $\varepsilon$
to be less than \foreignlanguage{ngerman}{$\frac{1}{2}$}, ie. $\varepsilon=\frac{3}{7}$.

In addition to the formula (64), a rough estimate of the order of
$\mathfrak{f}\left(s\right)$ for fixed $\sigma$ and $t\rightarrow\infty$
is still required. This may be obtained from the asymptotic development
of $\mathfrak{f}\left(s\right)$ in the region $t>0$, $-\sigma\leq t^{\frac{3}{7}}$.
A look at the proof of (64) shows, that one has applied to equation
(70) the condition $-\sigma\geq t^{\frac{3}{7}}$ only in the weaker
form $\sigma<\sigma_{0}$ where $\sigma_{0}$ is any real number.
It is valid therefore, in analogy to (70) with $r=0$, that
\begin{equation}
\mathfrak{f}\left(s\right)=\pi^{s-\frac{1}{2}}\dfrac{\varGamma\left(\dfrac{1-s}{2}\right)}{\varGamma\left(\dfrac{s}{2}\right)}\left(\zeta\left(1-s\right)-{\displaystyle \sum_{n=1}^{m}n^{s-1}}\right.
\end{equation}

\[
\left.-\eta^{s-1}e^{-\pi i\eta^{2}}\left\{ \dfrac{\sqrt{2}e^{\frac{3\pi i}{8}}\sin\pi\eta+\left(-1\right)^{m-1}e^{2\pi i\eta-2\pi i\left(\eta-m\right)^{2}}}{2\cos2\pi\eta}+O\left(\eta^{-1}\right)\right\} \right)
\]
with $\eta=\sqrt{\frac{s-1}{2\pi i}}$, $\left|\arg\eta\right|<\frac{\pi}{4}$,
$m=\left[\mathfrak{R\mathsf{\eta}+I}\mathsf{\eta}\right]$ in the
quarter plane $\sigma<\sigma_{0}$, $t>0$. Other members of the asymptotic
series may be obtained by the method of § 2, but are not required
for the present purpose.

A second asymptotic development of $\mathfrak{f}\left(s\right)$,
which is preferable for the quarter plane $\sigma>\sigma_{0}$, $t>0$,
one obtains by applying the saddle point method, not to the representation
of $\mathfrak{f}\left(s\right)$ provided by (65), but rather to that
of (58). The calculation need not be repeated in its entirety, since
the integral in (58) arises from that in (65), by passing to the complex
conjugate variables and by replacing $\sigma$ with $1-\sigma$. Consequently,
\begin{equation}
\mathfrak{f}\left(s\right)={\displaystyle \sum_{n=1}^{m_{1}}n^{s-1}}
\end{equation}
\[
+\eta_{1}^{-s}e^{\pi i\eta_{1}^{2}}\left\{ \dfrac{\sqrt{2}e^{-\frac{3\pi i}{8}}\sin\pi\eta_{1}+\left(-1\right)^{m_{1}-1}e^{-2\pi i\eta_{1}-2\pi i\left(\eta_{1}-m_{1}\right)^{2}}}{2\cos2\pi\eta_{1}}+O\left(\eta_{1}^{-1}\right)\right\} 
\]
\newpage\noindent with $\eta_{1}=\sqrt{\frac{s}{2\pi i}}$, $\left|\arg\eta_{1}\right|<\frac{\pi}{4}$,
$m_{1}=\left[\mathfrak{R\mathsf{\eta_{1}}-I}\mathsf{\eta_{1}}\right]$
in the quarter plane $\sigma>\sigma_{0}$, $t>0$. By comparing (82)
and (83) there results the asymptotic series for $\zeta\left(s\right)$
in each half-strip $\sigma_{1}<\sigma<\sigma_{2}$, $t>0$; this derivation
is perhaps somewhat simpler in relation to the necessary estimates
than that of § 2, but the individual terms of the series appear here
initially in a more complicated form.

From (83) it follows that
\begin{flalign}
 &  &  & \left\{ \begin{aligned}\mathfrak{f}\left(s\right) & ={\displaystyle \sum_{n=1}^{m_{1}}n^{-s}+O\left(\left(\dfrac{\left|s\right|}{2\pi e}\right)^{-\frac{\sigma}{2}}\right)} &  &  &  &  &  &  &  &  &  &  &  &  & \left(\sigma\geq0,t>0\right)\\
\mathfrak{f}\left(s\right) & =O\left(t^{\frac{1}{4}}\right) &  &  &  &  &  &  &  &  &  &  &  &  & \left(\sigma\geq\tfrac{1}{2}\right)\\
\left|\mathfrak{f}\left(s\right)\right.\! & \left.-\:1\right|<\dfrac{3}{4} &  &  &  &  &  &  &  &  &  &  &  &  & \left(\sigma\geq2,t>t_{0}\right).
\end{aligned}
\right.
\end{flalign}
and from (82)
\begin{flalign}
 &  &  & \left\{ \begin{aligned}\mathfrak{f}\left(s\right) & =\pi^{s-\frac{1}{2}}\,\dfrac{\varGamma\left(\dfrac{1-s}{2}\right)}{\varGamma\left(\dfrac{s}{2}\right)}\left(\zeta\left(1-s\right)-{\displaystyle \sum_{n=1}^{m}n^{s-1}+O\left(1\right)}\right) & \left(\sigma\leq1,t>0\right)\,\\
\mathfrak{f}\left(s\right) & =\pi^{s-\frac{1}{2}}\,\dfrac{\varGamma\left(\dfrac{1-s}{2}\right)}{\varGamma\left(\dfrac{s}{2}\right)}\;O\left(\left(\dfrac{t}{2\pi}\right)^{\frac{\sigma}{2}}\left|\sigma\right|^{-1}\right) & \left(0<-\sigma\leq t^{\frac{3}{7}},t>0\right)\\
\mathfrak{f}\left(s\right) & =\pi^{s-\frac{1}{2}}\,\dfrac{\varGamma\left(\dfrac{1-s}{2}\right)}{\varGamma\left(\dfrac{s}{2}\right)}\;O\left(\log t\right) & \left(0<-\sigma\leq t^{\frac{3}{7}},t>0\right).
\end{aligned}
\right.
\end{flalign}

In what follows it is convenient, instead of $\mathfrak{f}\left(s\right)$
to introduce the function
\[
g\left(s\right)=\pi^{-\frac{s+1}{2}}e^{-\frac{\pi is}{4}}\varGamma\left(\dfrac{s+1}{2}\right)\mathfrak{f}\left(s\right)
\]
By (64) one has, for $t>0$, $-\sigma\geq t^{\frac{3}{7}}$ and with
$\eta=\sqrt{\frac{s-1}{2\pi i}}$
\begin{flalign}
 &  &  &  &  & g\left(s\right)\sim e^{-\frac{7\pi i}{8}}\tan\dfrac{\pi s}{2}\dfrac{\sin\pi\eta}{\cos2\pi\eta}\qquad & \left(\left|s\right|\rightarrow\infty\right).
\end{flalign}
It should now be possible to estimate the mean value of $\left|g\left(s\right)\right|^{2}$
on each half-line $\sigma=\sigma_{0}<\frac{1}{2}$, $t\geq0$, specifically,
the expression
\begin{flalign*}
 &  &  &  &  & T^{-1}\intop_{0}^{T}\left|g\left(\sigma+ti\right)\right|^{2}dt\qquad\qquad & \left(\sigma<\tfrac{1}{2}\right).
\end{flalign*}
One may achieve this with the help of the asymptotic expansion (82);\newpage\noindent
the most elegant way is, however, a derivation using (58); from this,
specifically for \foreignlanguage{ngerman}{$\varepsilon>0$}
\[
\intop_{0}^{\infty}\left|\mathfrak{f}\left(\sigma+ti\right)\right|^{2}e^{-\varepsilon t}dt
\]
\[
=\intop_{0}^{\infty}e^{-\varepsilon t}\left\{ \,\;\intop_{\mathclap{0\swarrow1}}\;\;\;\intop_{\mathclap{0\searrow1}}\dfrac{x^{-\sigma-ti}y^{-\sigma+ti}e^{\pi i\left(x^{2}-y^{2}\right)}}{\left(e^{\pi ix}-e^{-\pi ix}\right)\left(e^{\pi iy}-e^{-\pi iy}\right)}dxdy\,\right\} dt,
\]
and here one may transform the right side by deformation of the path
of integration, interchange of the order of integration, and application
of the residue theorem. The calculation yields the statement
\[
\intop_{0}^{\infty}\left|\mathfrak{f}\left(\sigma+ti\right)\right|^{2}e^{-\varepsilon t}dt\sim\dfrac{1}{2\varepsilon}\left(2\pi\varepsilon\right)^{\sigma-\frac{1}{2}}\varGamma\left(\dfrac{1}{2}-\sigma\right),
\]
valid for $\sigma<\frac{1}{2}$ and $\varepsilon\rightarrow\infty$,
and from this follows further
\[
\intop_{1}^{\infty}\left|\mathfrak{f}\left(\sigma+ti\right)\right|^{2}\left(\dfrac{t}{2\pi}\right)^{\sigma}e^{-\varepsilon t}dt\sim\dfrac{\left(2\varepsilon\right)^{-\frac{3}{2}}}{1-2\sigma}.
\]
So for every fixed $\sigma<\frac{1}{2}$
\[
\intop_{1}^{T}\left|\mathfrak{f}\left(\sigma+ti\right)\right|^{2}\left(\dfrac{t}{2\pi}\right)^{\sigma}dt\sim\dfrac{1}{3\sqrt{2\pi}}\cdot\dfrac{T^{\frac{3}{2}}}{\frac{1}{2}-\sigma}.
\]
But from Stirling's formula comes the additional statement 
\begin{equation}
\left|g\left(s\right)\right|\sim\sqrt{2}\pi^{-\frac{\sigma}{2}}\left(\dfrac{t}{2}\right)^{\frac{\sigma}{2}}\left|\mathfrak{\mathfrak{f}}\left(s\right)\right|,
\end{equation}
and with this is won the desired formula
\[
T^{-1}\intop_{1}^{T}\left|g\left(\sigma+ti\right)\right|^{2}dt\sim\dfrac{1}{3}\sqrt{\dfrac{2}{\pi}}\cdot\dfrac{T^{\frac{1}{2}}}{\frac{1}{2}-\sigma}
\]
for fixed $\sigma<\frac{1}{2}$. From it follows further 
\begin{flalign}
 &  &  &  &  &  &  & \intop_{0}^{T}\log\left|g\left(\sigma+ti\right)\right|dt<\dfrac{T}{2}\log\dfrac{\sqrt{2}T^{\frac{1}{2}}}{3\sqrt{\pi}\left(\frac{1}{2}-\sigma\right)}+o\left(T\right) & \left(\sigma<\tfrac{1}{2},\:T\rightarrow\infty\right).
\end{flalign}

For $\sigma=\frac{1}{2}$ a lower bound for $\intop_{0}^{T}\log\left|g\left(\sigma+ti\right)\right|dt$
results. In fact, according to (60), on the critical line\newpage
\[
\left|\pi^{-\frac{s}{2}}\varGamma\left(\dfrac{s}{2}\right)\zeta\left(s\right)\right|\leq\left|2\pi^{-\frac{s}{2}}\varGamma\left(\dfrac{s}{2}\right)\mathfrak{f}\left(s\right)\right|,
\]
 thus, by (87)
\begin{flalign*}
 &  &  &  &  & \left|g\left(s\right)\right|\geq\left(8\pi\right)^{-\frac{1}{4}}t^{\frac{1}{4}}\left|\zeta\left(s\right)\right|\left(1+o\left(1\right)\right)\qquad & \left(\sigma=\tfrac{1}{2}\right)
\end{flalign*}
\begin{equation}
\intop_{0}^{T}\log\left|g\left(\frac{1}{2}+ti\right)\right|dt>\dfrac{T}{4}\log T-\left(\log8\pi+1\right)\dfrac{T}{4}
\end{equation}
\[
+\intop_{0}^{T}\log\left|\zeta\left(\frac{1}{2}+ti\right)\right|dt+o\left(T\right).
\]

Finally, for $\sigma\geq2$, by (87) and (84)
\begin{equation}
\intop_{0}^{T}\log\left|g\left(\sigma+ti\right)\right|dt=\sigma\left(\dfrac{T}{2}\log\dfrac{T}{2\pi}-\dfrac{T}{2}\right)+\dfrac{T}{2}\log2+o\left(T\right).
\end{equation}

Now let $t_{0}>0$, $T>t_{0}$, and the straight line $t=t_{0}$,
$t=T$   be free of zeros of the function $g\left(s\right)$. Furthermore,
let $\sigma_{0}>-T^{\frac{3}{7}}=\sigma_{1}$. Consider the rectangle
with sides $\sigma=\sigma_{0}$, $t=T$, $\sigma=\sigma_{1}$, $t=t_{0}$.
On the left side $\sigma=\sigma_{1}$, $t=t_{0}\leq t\leq T$ there
are, for sufficiently large $T$, by (64) no zeros of $g\left(s\right)$.
One connects the zeros of $g\left(s\right)$ lying within the rectangle
through cuts that are constructed parallel to the real axis, with
the right side $\sigma=\sigma_{0}$. In the cut rectangle $\log\,g\left(s\right)$
is then unambiguous {[}\foreignlanguage{ngerman}{eindeutig}{]}; it
becomes a branch of this function by requiring that $0\leq\arg\,g\left(\sigma_{1}+Ti\right)<2\pi$
is fixed. As is known, then it applies that 
\begin{align}
2\pi\sum_{\alpha<\sigma_{0}}\left(\sigma_{0}-\alpha\right) & =\intop_{t_{0}}^{T}\log\left|g\left(\sigma_{0}+ti\right)\right|dt-\intop_{\sigma_{1}}^{\sigma_{0}}\arg\,g\left(\sigma+Ti\right)d\sigma\\
 & -\intop_{t_{0}}^{T}\log\left|g\left(\sigma_{1}+ti\right)\right|dt+\intop_{\sigma_{1}}^{\sigma_{0}}\arg\,g\left(\sigma+t_{0}i\right)d\sigma,\nonumber 
\end{align}
where $\alpha$ passes through the real parts of all roots of $g\left(s\right)$
located in the square. The first integral can be accurately estimated
for $\sigma_{0}<\dfrac{1}{2}$ upward, for $\sigma_{0}=\dfrac{1}{2}$
downward and for $\sigma_{0}\geq2$. The third and fourth integral
contribute, as can be seen without significant difficulty by (86),
only an amount of the order $T^{\frac{13}{14}}$. Finally, the second
integral can be estimated in the usual manner by using (84) and (85)
as $O\left(T^{\frac{6}{7}}\log T\right)$. Accordingly, by (88) \newpage
\begin{flalign}
 &  &  &  &  & \sum_{\alpha<\sigma}\left(\sigma-\alpha\right)<\dfrac{T}{8\pi}\log T-\dfrac{T}{4\pi}\log\left\{ 3\sqrt{\dfrac{\pi}{2}}\left(\dfrac{1}{2}-\sigma\right)\right\} +o\left(T\right) & \left(\sigma<\tfrac{1}{2}\right),
\end{flalign}
 by (89) 
\begin{equation}
\sum_{\alpha<\frac{1}{2}}\left(\dfrac{1}{2}-\alpha\right)<\dfrac{T}{8\pi}\log T-\left(1+\log8\pi\right)\dfrac{T}{8\pi}+\dfrac{1}{2\pi}\intop_{0}^{T}\log\left|\zeta\left(\frac{1}{2}+ti\right)\right|dt+o\left(T\right)
\end{equation}
and by (90)
\begin{flalign}
 &  &  &  &  & \sum_{\alpha}\left(\sigma-\alpha\right)=\sigma\left(\dfrac{T}{4\pi}\log\dfrac{T}{2\pi}-\dfrac{T}{4\pi}\right)-\dfrac{T}{4\pi}\log2+o\left(T\right) & \left(\sigma\geq2\right);
\end{flalign}
in the last equation by running $\alpha$ through the real parts of
all zeros of $g\left(s\right)$ lying in the strip $0<t<T$. If their
number is designated by $N_{1}\left(T\right)$, then by (94)
\begin{equation}
N_{1}\left(T\right)=\dfrac{T}{4\pi}\log\dfrac{T}{2\pi}-\dfrac{T}{4\pi}+o\left(T\right).
\end{equation}
In the upper half plane, the zeros of $g\left(s\right)$ agree with
those of $\mathfrak{f}\left(s\right)$. Under the condition that,
of the $N_{1}\left(T\right)$ zeros, up to $o\left(T\right)$ of them
lie to the right of $\sigma=\frac{1}{2}$, then the change of $\arg\,\mathfrak{f}\left(\frac{1}{2}+ti\right)$
in the interval $0<t<T$ would equal $-\left(\dfrac{T}{2}\log\dfrac{T}{2\pi}-\dfrac{T}{2}\right)+o\left(T\right)$,
and one does not obtain a statement about the zeros of $\zeta\left(\frac{1}{2}+ti\right)$.
First, however, it follows from (94), further, that
\[
\sum_{\alpha}\alpha=-\dfrac{T}{4\pi}\log2+o\left(T\right);
\]
so there are certainly an infinite number of zeros of $\mathfrak{\mathfrak{f}}\left(s\right)$
that even lie to the left of $\sigma=0$; and this is obtained from
(92) and (93), independently of (94), by subtracting a lower bound
for the number of zeros of $\mathfrak{\mathfrak{f}}\left(s\right)$
lying in the region $\sigma<\frac{1}{2}$, $0<t<T$. We denote this
number by $N_{2}\left(T\right)$, so it follows for each $\sigma<\frac{1}{2}$
\[
\left(\dfrac{1}{2}-\sigma\right)N_{2}\left(T\right)>\dfrac{T}{4\pi}\log\left\{ \dfrac{3}{4}e^{-\frac{1}{2}}\left(\dfrac{1}{2}-\sigma\right)\right\} +\dfrac{1}{2\pi}\intop_{0}^{T}\log\left|\zeta\left(\frac{1}{2}+ti\right)\right|dt+o\left(T\right).
\]
This estimate is most favorable for
\[
\sigma=\dfrac{1}{2}-\dfrac{4}{3}e^{\frac{3}{2}}
\]
and results in
\[
N_{2}\left(T\right)>\dfrac{3}{16\pi}e^{-\frac{3}{2}}T+\dfrac{3}{8\pi}e^{-\frac{3}{2}}\intop_{0}^{T}\log\left|\zeta\left(\frac{1}{2}+ti\right)\right|dt+o\left(T\right).
\]
\newpage{}

It now is well-known that, using an approach of the form of (91),
with $\zeta\left(s\right)$ in place of of $g\left(s\right)$, 
\[
\frac{1}{2\pi}\intop_{0}^{T}\log\left|\zeta\left(\frac{1}{2}+ti\right)\right|dt={\displaystyle \sum_{\alpha_{\zeta}>\frac{1}{2}}}\left(\alpha_{\zeta}-\frac{1}{2}\right)+O\left(\log T\right),
\]
where $\alpha_{\zeta}$ passes through the real parts of the zeros
of the zeta function lying in the strip $0<t<T$ located to the right
of the critical line. It follows that
\begin{equation}
N_{2}\left(T\right)>\dfrac{3}{16\pi}e^{-\frac{3}{2}}T+\dfrac{3}{4}e^{-\frac{3}{2}}{\displaystyle \sum_{\alpha_{\zeta}>\frac{1}{2}}}\left(\alpha_{\zeta}-\frac{1}{2}\right)+o\left(\log T\right).
\end{equation}
To the right of $\sigma=\frac{1}{2}$ there are at most $N_{1}\left(T\right)-N_{2}\left(T\right)$
zeros of $\mathfrak{f}\left(s\right)$ within the strip $0<t<T$,
so $\arg\,\mathfrak{f}\left(\frac{1}{2}+ti\right)$ is reduced in
the interval $0<t<T$ by at most $2\pi\left(N_{1}\left(T\right)-N_{2}\left(T\right)\right)+O\log\left(T\right)$.
Consequently, $\arg\,\varphi\left(\frac{1}{2}+ti\right)$ increases
in this interval by at least

\[
\vartheta\left(T\right)-2\pi N_{1}\left(T\right)+2\pi N_{2}\left(T\right)+O\left(\log T\right),
\]
and this quantity is by (63), (95), (96) at least equal to $2\pi N_{2}\left(T\right)+o\left(T\right)$.
Therefore $N_{0}\left(T\right)$, the number of zeros of $\zeta\left(\frac{1}{2}+ti\right)$
in the interval $0<t<T$, satisfies the inequality
\begin{equation}
N_{0}\left(T\right)>\dfrac{3}{8\pi}e^{-\frac{3}{2}}T+\dfrac{3}{2}e^{-\frac{3}{2}}{\displaystyle \sum_{\alpha_{\zeta}>\frac{1}{2}}}\left(\alpha_{\zeta}-\frac{1}{2}\right)+o\left(T\right).
\end{equation}
The density of zeros of $\zeta\left(s\right)$ lying on the critical
line, that is, the lower bound of the ratio $N_{0}\left(T\right):T$
for $T\rightarrow\infty$, therefore is positive and at least equal
to $\dfrac{3}{8\pi}e^{-\frac{3}{2}}$, thus is greater than $\frac{1}{38}$.
Aside from this numerical value, the statement is not new, but has
already been demonstrated in 1920 by Hardy and Littlewood in a much
simpler way. Nonetheless, this minor result may perhaps suggest some
independent value for the properties of $\mathfrak{f}\left(s\right)$
discussed here.

One may make an additional remark regarding formula (97). The falsehood
of Riemann's conjecture may be measured, in a certain sense, by the
sum $\sum\left(\alpha_{\zeta}-\dfrac{1}{2}\right)$. Although one
knows, through Littlewood, that this total does not exceed $O\left(T\log\log T\right)$,
no better estimate is known. If the Riemann hypothesis is false, this
sum may grow faster than $T$; then, however, by (97), the number
$N_{0}\left(T\right)$ would grow faster than $T$, and the Riemann\newpage\noindent
hypothesis could not be ``too false''. Let $\psi\left(t\right)$
stand for any positive function of $t$ which grows toward infinity
more slowly than $\log t$, then it still follows from (97) that,
in the narrow region $0\leq\sigma-\frac{1}{2}\leq\frac{\psi\left(t\right)}{\log t}$,
$2\leq t\leq T$, there are at least $\frac{3}{4\pi}e^{-\frac{3}{2}}T\psi\left(T\right)\left(1+o\left(1\right)\right)$
zeros of $\zeta\left(s\right)$. This is a new result, valid even
for the case where $\psi\left(t\right)$ increases more slowly than
$\log\log t$. So there are, for example, in the region $0\leq\sigma-\frac{1}{2}\leq\frac{19}{\log t}$
, $2\leq t\leq T$ more than $T+o\left(T\right)$ zeros.

The question remains open as to whether one can improve the lower
bound of $N_{2}\left(T\right)$ given in (96). For the proof of Riemann's
assertion that $N_{0}\left(T\right)$ is asymptotically equal to $\frac{T}{2\pi}\log\frac{T}{2\pi}-\frac{T}{2\pi}$,
it suffices to show the corresponding statement for $N_{2}\left(T\right)$.
This seems difficult to achieve with the previously used methods of
analysis in the theory of the zeta function without an essentially
new idea; and this is especially true of any attempts to prove the
Riemann Hypothesis.
\end{document}